\documentclass[reqno,centertags,12pt]{amsart}

\usepackage{amssymb,amsmath,amsthm,dsfont}
\usepackage{graphicx}
\usepackage{color}
\usepackage{bigstrut}
\makeatletter \def\@strippedMR{} \def\@scanforMR#1#2#3\endscan{%
  \ifx#1M\ifx#2R\def\@strippedMR{#3}%
  \else\def\@strippedMR{#1#2#3}%
  \fi\fi} \renewcommand\MR[1]{\relax \ifhmode\unskip\spacefactor3000
  \space\fi \@scanforMR#1\endscan
  MR\MRhref{\@strippedMR}{\@strippedMR}} \makeatother
\newcommand{\R}{\mathbb{R}} \newcommand{\Z}{\mathbb{Z}}
\newcommand{\T}{\mathbb{T}} \newcommand{\C}{\mathbb{C}}

\renewcommand{\^}{\widehat}

\theoremstyle{plain}
\newtheorem{thm}{Theorem}[section]
\newtheorem{cor}[thm]{Corollary}
\newtheorem{lem}[thm]{Lemma}

\newtheorem{defin}[thm]{Definition}

\theoremstyle{definition}

\newtheorem{rem}{Remark}

\newcommand{\beq}{\begin{equation}}
\newcommand{\eeq}{\end{equation}}

\def\pasdegrille{\let\grille = \pasgrille}
\def\ecriture#1#2{\setbox1=\hbox{#1}
\dimen1= \wd1 \dimen2=\ht1 \dimen3=\dp1 \grille #2 \box1 }
\def\aat#1#2#3{
\divide \dimen1 by 48 \dimen3=\dimen1 \multiply \dimen1 by #1
\advance \dimen1 by -\dimen3 \divide \dimen1 by 101 \multiply
\dimen1 by 100 \divide \dimen2 by \count11 \multiply \dimen2 by #2
\setbox0=\hbox{#3}\ht0=0pt\dp0=0pt
  \rlap{\kern\dimen1 \vbox to0pt{\kern-\dimen2\box0\vss}}\dimen1= \wd1
\dimen2=\ht1}
\def\pasgrille{
\count12= \dimen1 \divide \count12 by 50 \divide \dimen2 by \count12
\count11 =\dimen2 \ \divide \dimen1 by 48
\setlength{\unitlength}{\dimen1} \smash{\rlap{\ }} \dimen1= \wd1
\dimen2=\ht1 }
\def\grille{
\count12= \dimen1 \divide \count12 by 50 \divide \dimen2 by \count12
\count11 =\dimen2 \ \divide \dimen1 by 48
\setlength{\unitlength}{\dimen1}
\smash{\rlap{\graphpaper[1](0,0)(50, \count11)}} \dimen1= \wd1
\dimen2=\ht1 }

\pasdegrille

\def\a{\alpha}
\def\rta{\rightarrow}
\def\q{\quad}
\def\e{\varepsilon}
\def\t{\tau}
\def\d{\delta}
\def\sg{\sigma}
\def\ld{\lambda}
\def\z{\zeta}
\def\C{\mathcal{C}}
\def\G{G_{\t_1,\,\t_2,\,\t_3}}
\def\T{\mathcal{T}}

\begin{document}

\title[2D Maximal inequality]
{An improved maximal inequality for 2D  fractional order
Schr\"{o}dinger operators}

\author[Miao et al]{Changxing Miao, Jianwei Yang and Jiqiang Zheng}


\address{Institute of Applied Physics and Computational Mathematics, Beijing 100088, China,
\;\;
Beijing Center of Mathematics and Information Sciences,
Beijing, 100048, P.R.China.}
\email{miao\_{}changxing@iapcm.ac.cn}

\address{Beijing International Center for Mathematical Research, Peking University, Beijing, 100871, P.R.China} \email{geewey\_young@pku.edu.cn}

\address{Universit\'e Nice Sophia-Antipolis, 06108 Nice Cedex 02, France}
\email{zhengjiqiang@gmail.com}

\begin{abstract}
  The local maximal inequality for the
  Schr\"{o}dinger operators of order $\a>1$ is shown to be bounded
  from $H^s(\R^2)$ to $L^2$ for any $s>\frac38$.
  This improves the previous result of Sj\"{o}lin
  on the regularity of solutions to fractional order Schr\"{o}dinger equations.
  Our method is inspired by Bourgain's argument in case of $\a=2$.
  The extension from $\a=2$ to general $\a>1$ confronts three
  essential obstacles: the lack of Lee's reduction lemma,
  the absence of the algebraic structure of the symbol
  and the inapplicable Galilean transformation in the deduction of the main theorem.
  We get around these difficulties by establishing a new reduction lemma at our disposal
  and analyzing all the possibilities in using the separateness of the segments
  to obtain the analogous bilinear $L^2-$estimates.
  To compensate the absence of Galilean invariance,
  we resort to Taylor's expansion for the phase function.
  The Bourgain-Guth inequality in \cite{ref Bourgain Guth} is also rebuilt
  to dominate the solution of fractional order Schr\"{o}dinger equations.
\end{abstract}
\maketitle
\tableofcontents
\section{Introduction and the main result}\label{sect:intro_main}
\noindent

For $\a>1$, we define the $\a$-th order Schr\"{o}dinger evolution
operator by
$$
 U(t)f(x)
 \mathrel{\mathop=^{\rm def}}
 (2\pi)^{-\frac d 2}\int_{\R^d}e^{i[\,x\cdot\xi+t\,|\xi|^\a]}\^{f}(\xi)\,d\xi,
$$
and consider the following local maximal inequality
\begin{equation}\label{1}
 \Bigl\|\sup_{0<t<1}|U(t)f|\Bigr\|_{L^2(B(0,1))}\leq C_{\a,s,d}\|f\|_{H^{s}(\R^d)},\;s\in
 \R\,,
\end{equation}
where $H^s(\R^d)$ is the usual inhomogeneous Sobolev space defined
via Fourier transform, and $B(0,1)$ is the unit ball
centered at the origin. As a consequence of \eqref{1}, we have
the following point-wise convergence phenomenon from a standard process of approximation
and Fatou's lemma
$$
 \lim_{t\rta 0}U(t)f(x)=f(x),\q\text{a.e.}\;x,\q \forall\, f\in H^s(\R^d).
$$
If $s>\frac{d}{2}$, we obtain \eqref{1} immediately from Sobolev's
embedding. Thus, it is natural to ask what the minimal $s$ is
to ensure $\eqref{1}$.\vskip0.2cm

First of all, let us briefly review the known results about \eqref{1}
in the case when $\a=2$. This problem was raised by Carleson in \cite{ref Car}, where
he answered the 1D case with $s\geq\frac14$. This result was shown
to be optimal by Dahlberg and Kenig~\cite{ref DK}. The higher dimensional cases of
\eqref{1} were established independently by Sj\"{o}lin \cite{ref Sj}
and Vega \cite{ref Vega} with $s>\frac12$. In particular, the result
can be strengthened to $s=\frac12$ by Sj\"{o}lin \cite{ref Sj} when
$d=2$. Meanwhile, Vega \cite{ref Vega} demonstrated that \eqref{1}
fails in any dimension if $s<\frac14$. It is then conjectured that
$s\geq\frac14$ should be sufficient for all dimensions.

The breakthrough  was achieved by Bourgain \cite{ref Bourgain1,ref
Bourgain2}, where he showed that \eqref{1} holds with $\a=2$ for some
$s<\frac12$ when $d=2$. His work was carried over and improved
subsequently by many authors, including  Moyua, Vargas, Vega
\cite{ref M-V-V}, Tao and Vargas \cite{ref T-V1,ref T-V2}, Lee
\cite{ref Lee} and Shao \cite{ref shao}, where the best result
hitherto is $s>\frac{3}{8}$ due to Lee \cite{ref Lee}.

Previous to \cite{ref Bourgain3}, the results about $d\geq 3$ remained
$s>\frac12$, and $s\geq\frac14$ was still believed to be the correct
condition for \eqref{1} in every dimension. The study on this
problem stagnated for several years until the recent work \cite{ref
Bourgain3}, where the $\frac12-$barrier was broken for \emph{all}
dimensions. More surprisingly, Bourgain also discovered some
counterexamples to disprove the widely believed assertion on the
$\frac14-$threshold. Specifically, he showed that
$s\geq\frac12-\frac1d$ is necessary for \eqref{1} if $d\geq 4$.
These examples originated essentially from an
observation on arithmetical progressions.\vskip0.2cm

Now, let us come to the fractional order case. Sj\"{o}lin proved
$\eqref{1}$ with $s=\frac12$, $d=2$ for \emph{all} $\a>1$ in \cite{ref Sj}.
His proof involves a $TT^*$ argument, which reduces the
problem to a dispersive estimate of a specific oscillatory
integral. After localizing the integral to the high and low
frequencies, the author employed a classical result by Miyachi \cite{ref
Miyachi} to treat the high frequency part. The other part is
estimated by means of the following inequality $$
\int_{\R^2}|y|^{-1}(1+|x-y|^4)^{-1}dy\leq\, c\, |x|^{-1},
$$
where the decay rate in the right hand side can not be
improved.
A crucial fact which Sj\"{o}lin's proof relied heavily on is
the factor $|t(x)-t(y)|^{\frac1\a}$ can be canceled at the end of
the computation exactly for $s=\frac12$. This is usually
referred as the Kolmogoroff-Seliverstoff-Plessner method, see
\cite{ref Car} and \cite{ref Sj} for more details.
Due to these reasons, it seems very difficult to
pursue Sj\"{o}lin's original approach to improve this
result. In this paper, we prove
\begin{thm}\label{thm 1}
 If $d= 2$ and $\a>1$, then \eqref{1} is valid for all $s>\frac{3}{8}$.
\end{thm}

\begin{rem}\label{rmk 1}
As noted in \cite{ref Bourgain3}, one may modify the method to treat
a general multiplier operator $\Phi(D)$ having the property that
for some constants $C,\,c>0$ and all multi-indices
$\gamma$
$$
 |\partial^\gamma\Phi(\xi)|\leq\,C|\xi|^{2-|\gamma|},\q|\nabla\Phi(\xi)|\geq c\,|\xi|.
$$
However, this does not concern the fractional order case.
\end{rem}

As a consequence, we get some improvement on the higher dimensional results
by using the scheme of induction on dimensions
formulated in \cite{ref Bourgain3}.
\begin{cor}\label{corol}
 For $d\geq 3$ and $\a>1$, there exists a $\theta_d$ such that \eqref{1} is valid for all
 $s>\theta_d$ with
 $$
  \theta_d=\frac{1}2-\sigma\Big(\frac{1}{2}-\theta_{d-1}\Big),
 $$
 for some $\sigma\in (0,\frac12)$.
 In particular, $\theta_d<\frac{1}{2}$ for every $d\geq 2$ since $\theta_2<\frac12$.
\end{cor}
\begin{rem}\label{rmk 2}
This improves Theorem 2 in \cite{ref Sj} in higher dimensions. Noting that
the induction argument in \cite{ref Bourgain3} is independent of the
order $\a$, we may apply it verbatim to obtain Corollary
\ref{corol}.
\end{rem}

As in \cite{ref Bourgain3}, the proof is based on the multilinear restriction theorem in \cite{ref BCT}.
To achieve this, an important observation introduced by Bourgain and Guth \cite{ref Bourgain Guth} is that
up to an $R^\e$ factor and  a well behaved remainder, one successfully control the free solution of the Schr\"{o}dinger equation with
a summation of triple products fulfilling the transversality condition for which the multilinear restriction estimate
can be used. Roughly speaking, one gains structures by losing $R^\e$, however this is acceptable if we do not intend to
solve the end-point problem.
These triple products, which we will call type I terms in Section 5,
are generated by iteration with respect to scales.
As a result, they are used to collect the contributions obtained at different scales.
In this sense, it is also reminiscent of Wolff's induction on scale argument in \cite{ref Wolff}.\footnote{It is thus interesting to consider how to combine these two
important ideas together to improve the argument in this work.}
In this paper, we call this robust device as Bourgain-Guth's inequality.
\vskip0.2cm

Let us take this opportunity to try to moderately
clarify several points in Bourgain's argument reserving the notations in \cite{ref Bourgain Guth} and \cite{ref Bourgain3}. Of course, we do not have the ambition to present a complete clarification of Bourgain's treatment since that will be far beyond our reach. Instead, we focus only on the points which are directly relevant to this paper. First, a crucial input is the Bourgain-Guth inequality for oscillatory integrals from \cite{ref Bourgain Guth}. It collects the contributions of the transversal triple products from all dyadic scales between $R^{-\frac12}$ and one, so that we can use the multilinear restriction theorem in \cite{ref BCT} to evaluate the contributions at each scale. Since we are dealing with dyadic scales in $(R^{-\frac12},1)$, we can safely consider items from all scales by taking an $\ell^2$ sum losing at most a factor $\log R$.
To obtain this inequality, they tactically used a `` local constant trick" in \cite{ref Bourgain Guth} according to the following principle.  By writing the oscillatory integral into trigonometric sums with variable coefficients $T_\a f(x)$, one may regard $T_\a f(x)$ as a constant on each ball of radius $K$ thanks to the uncertainty principle, where these ``constants" certainly depend on the position of the ball. This  heuristic point is justified by convolving $T_\a f$ with some suitable bump functions. However, to carry out further manipulations especially the iterative process, it is awkward to write it out explicitly and repeatedly. Instead of this, one prefers doing formal calculation for brevity and clarity. Based on this observation, one may insert/extract the factor $T_\a f(x)$ into/from an integral over a ball of radius $K$, or more generally over a tile in suitable shape and size. All the judgements are made according to the uncertainty principle.
This simple and important observation is very efficient in simplifying various explicit calculations in the  context so that the Bourgain-Guth inequality can be established in \cite{ref Bourgain Guth} by iteration. Let us say more about the establishment of Bourgain-Guth's inequality before turning to the argument for  the Schr\"odinger maximal function. The brilliant novelty in \cite{ref Bourgain Guth}, which we will follow in Section \ref{sect:key decomp}, embodied in the way of using Bonnet-Carbery-Tao's multiplier restriction theorem. The idea might be roughly described as after writing $Tf(x)$ into a variable coefficient trigonometric sum, one may estimate for each $x\in B_R$ in three different manners, where only a small portion of the members in $\{T_\a f(x)\}_\a$ would dominate the behavior of $Tf(x)$. As can be seen in \cite{ref Bourgain Guth} and Section \ref{sect:key decomp}, these members correspond  respectively to three different scenarios which covers all the possibilities for a particular $x\in B_R(0)$ to encounter. According to \cite{ref Bourgain Guth}, they are titled as \emph{ non-coplanar interaction, non-transverse interaction} and \emph{transverse coplanar interaction}. We refer to Section \ref{sect:key decomp} for more detailed discussions about this classifications.
Now we turn to Bourgain's treatment on Schr\"odinger maximal function. The idea is that by using Bourgain-Guth's inequality, one is reduced to controlling each item in the $\ell^2$ summation with desired bound. To achieve this, one tiles $\R^3$ with translates of polar sets of the cap $\tau$, which contains a triple of transversal subcaps $\tau_1,\tau_2,\tau_3$. This provides a decomposition of $B_R\subset\R^3$. Invoking the local constant principle, one may raise and lower the moment exponents on each tile so that the favorite trilinear restriction in \cite{ref BCT} can be used. During this calculation, Galilean's transform is employed to shift the center of the square where the frequency is localized to the center. Although we have compressed Bourgain's argument into as few words as possible, it is far more difficult and subtle in concrete manipulations as in \cite{ref Bourgain3,ref Bourgain Guth}. We confine ourselves with this brief investigation on Bourgain's approach and turn to our situation below.
\vskip0.2cm

 To use the strategy in \cite{ref Bourgain3},
we need rebuild the Bourgain-Guth inequality for general $\a>1$.
Although this inequality is invented in \cite{ref Bourgain Guth} for $\a=2$,
it is rather non-trivial to generalize this result to $\a>1$ as will be seen in Section 5.
One of the obstructions is the absence of the algebraic structure of $|\xi|^\a$ when $\a$ is not an integer.
This fact leads to the distinctions of our argument from \cite{ref Bourgain3} and \cite{ref Bourgain Guth} in almost every aspects,
especially in the proof of the bilinear $L^2-$estiamte in Subsection 5.2 where we introduce a new argument.
\vskip0.2cm

Besides the reestablishment of Bourgain-Guth's inequality,
we need a fractional order version of Lee's reduction lemma  in \cite{ref Lee}
for general  $\a>1$. In Section 4, we establish this result using a different method.
This extends the result in  \cite{ref Lee} to a more general setting.
We will use the method of stationary phase in spirit of \cite{ref shao}.
However, to justify the proof, we involve a localization argument
which eliminates Schwartz tails by losing $R^\e$.
To be precise, we separate the Poisson summation
to a relatively large and small scales, where either the rapid decreasing property of Schwartz functions
or the stationary phase argument can be used to handle the error terms.
This principle is also used in the proof of
the main theorem. The essence of this argument is exploiting the orthogonality in ``phase space" via stationary phase
and Poisson summation formula. In doing this, one only needs to afford an $R^\e$ loss but one may sum the pieces that
are well-estimated efficiently. See Section 3 and Section 4 for more details.
\vskip0.2cm

At the end of this section, let us say a word about the potential of Bourgain-Guth's approach to oscillatory integrals.
In doing harmonic analysis, one of the most important principles is that
structures are favorable conditions to help us use deep results in mathematics.
For instance, Whitney's decomposition was employed by the authors in \cite{ref T-V1,ref T-V2,ref T-V-V 1998} to generate the transversality conditions
for the use of bilinear estimates. On the other hand, the proof of
Bourgain and Guth's inequality enlightened a new approach
to generate structures by means of logical classification, \emph{i.e.}
exploiting the intrinsic structures implicitly involved in the the
summation of large number of elements creatively using logical division.
The idea is fairly new and the argument is really a \emph{tour de force},
bringing in ideas and techniques from combinatorics as will be seen in Section 5.
We believe this approach is very promising to get improvements on the open questions in classical harmonic analysis.
In particular, the result in this paper might be improved further by refining this method.
\vskip0.2cm

This paper is organized as follows. In Section 2, we introduce some preliminaries and
the basic lemmas. In Section 3, we prove the main result. Section 4 is devoted to the proof of Lemma \ref{lem 1}
and Section 5 is devoted to the proof of Lemma \ref{key_decomp}.

\section{Preliminaries}\label{sect:not_ls}
\noindent

This section includes the list of the frequently used notations,
the statement of the crucial lemma which plays the
key role in deducing the main result as well as
the primary reduction for the proof of Theorem \ref{thm 1}.

\subsection{Notations}
\noindent

Throughout this paper the following notations will be used.
\\
$\diamond$ $\Omega=[-1/2,1/2]\times[-1/2,1/2]$
\\
$\diamond$ $[r]$ is the greatest integer not exceeding $r$.
\\
$\diamond$ If $\Omega$ is a subset of $\R^d$, we define $\displaystyle
\Omega^c=\R^d\setminus \Omega$.
\\
$\diamond$ $\chi_{\Omega}$ denotes the characteristic
function of a set $\Omega\subset\R^d$.
\\
$\diamond$ Suppose $\xi$ is a vector in $\R^d$, we
define $\overline{\xi}=\frac{\xi}{|\xi|}$.
\\
$\diamond$ We define $\mathcal{I}=\{\xi\in \R^d\mid 1/2\leq |\xi|\leq 2\}$ and,
except Lemma \ref{lem 1} below, we always assume $d=2$.
\\
$\diamond$ For $f(x)$ a measurable function and $a\in\R^d$, we define
$\displaystyle\t_a f(x)=f(x-a)$.
\\
$\diamond$  Denote by $ \mathcal{S}(\R^d)$ the Schwartz class on $\R^d$
and by $ \mathcal{S}'(\R^d)$ the space of tempered distributions.
\\
$\diamond$ We use $\mathcal{F}_{x\rta \xi}f$ or $\widehat{f}(\xi)$
to denote the Fourier transform of a tempered distribution $f(x)$ and
$$
\^f(\xi)=(2\pi)^{-\frac d2}\int_{\R^d}f(x)e^{-ix\cdot\xi}dx.
$$
\\
$\diamond$ Denote by $B(a,K)$ or $B_{a, K}$ a ball centered at $a$
in $\R^d$ of radius $K$.
\\
$\diamond$ Suppose $B$ is a convex body in $\R^d$ and $\ld>0$,
we use $\ld B$ to denote the convex set having the same center with $B$ but
enlarged in size by $\ld.$
\\
$\diamond $ The capital letter $C$ stands for a constant which
might be different from line to line and $c\ll C$ means $c$ is far less than
$C$. This is clear in the context.
\\
$\diamond $ The notion $A\lesssim B$ means $A\leq CB$ for some
constant $C$, and $A\simeq B$ means both $A\lesssim B$ and
$B\lesssim A$.
\\
$\diamond$ By $A\lesssim_{\eta,\zeta,\ldots}B$, we mean there is a
constant $C=C(\eta,\zeta,\ldots)$ depending on $\eta,\zeta,...$ such
that $A\leq C B$, where the reliance of $C$ on the parameters will
be clear from the context.
\\
$\diamond$ By $\zeta=\mathcal{O}_\a(\eta)$, we mean there is an estimate
$\zeta\lesssim_\a\eta$.

\subsection{Caps, tiles and the Bourgain-Guth inequality}
\noindent

Now we introduce some terminologies. Let $R\gg 2^{5\a}> 1$ and
$\frac{1}{\sqrt{R}}<\d<1$. Partition $\R^2$ into $\cup_\t\Omega_\t$
where $\Omega_\tau$ is a $\d\times\d$ square centered at $\xi_\tau
\in\d\Z^2$ such that the edges of $\Omega_\t$ are parallel to the
abscissas and vertical axis respectively. Let
$\overrightarrow{n}_\tau$ be the exterior unit normal vector of the
immersed surface $(\xi,\,|\xi|^\a)$ at the point
$(\xi_\t,\,|\xi_\tau|^\a)$. We define the following sets
\begin{align*}
\Pi^\d_\t&=(\xi_\t,\,|\xi_\t|^\a)+\left\{z\in\R^3\mid
|\langle z,\overrightarrow{n}_\t\rangle|\leq \d^\a\right\},\\
\mathcal{C}_\t&=\Pi^\d_\t\cap\left\{z\in \R^3\mid\, z=(z_1,z_2,z_3),\,
(z_1,\,z_2)\in\Omega_\t\right\}.
\end{align*}
Obviously,\, $\mathcal{C}_\t$ is a parallelopiped with dimensions
$\sim \d,\,\d,\,\d^\a$.

\begin{defin}
The parallelopiped $\mathcal{C}_\t$ is called a $\d-$cap
associated to $\Omega_\t$.
\end{defin}
\begin{defin}
The polar set of $\mathcal{C}_\t$ is defined as
$$
\mathcal{C}^*_\t=\Bigl\{z\in\R^3\;\Big|\; |\langle
z,\,w\rangle|\leq 1,\,\forall\, w\in\mathcal{C}_\t-(\xi_\t,\,|\xi_\t|^\a)\Bigr\},
$$
\end{defin}
It is easy to see that $\mathcal{C}^*_\t$ is essentially a
$\frac{1}{\d}\times\frac{1}{\d}\times\frac{1}{\d^\a}$-rectangle
centered at the origin, with the longest side in the direction of
$\overrightarrow{n}_\t$. Moreover, we may tile $\R^3$ with boxes of the
translations  of $\mathcal{C}^*_\t$. This decomposes
$\R^3$ naturally into the union of essentially disjoint $\mathcal{C}^*_\t-$boxes.
We call this decomposition a tiling of $\R^3$ with $\mathcal{C}^*_\t-$boxes.
\vskip0.2cm
Define an oscillatory integral by
$$
Tf(x)=\int_{\mathcal{I}}
e^{i[x_1\xi_1+x_2\xi_2+x_3|\xi|^\a]}\widehat{f}(\xi)\,d\xi,
$$
where $x=(x_1,\,x_2,\,x_3)\in\R^3$ and $\xi=(\xi_1,\xi_2)\in\R^2$.
Setting $x'=(x_1,x_2)$ and
regarding $x_3$ as the temporal
variable $t=x_3$, we have
\begin{equation}\label{osc-sol}
U(t)f(x')=Tf(x).
\end{equation}
\begin{rem}
The notion in \eqref{osc-sol} for general $d-$dimensional counterpart is defined
in the same way, with $x_{d+1}$ in place of $x_3$, and
$(x_1,x_2)$ as well as $(\xi_1,\xi_2)$ replaced by $(x_1,\ldots,x_{d})$ and $(\xi_1,\ldots,\xi_{d})$.
This will only be used for Lemma \ref{lem 1} below, which is proved for general dimensions.
\end{rem}

Now we can state Bourgian-Guth's inequality which will be used to control the
oscillatory integral $Tf(x)$ in terms of $\{Tf_\t\}_\t$, where
$\widehat{f}_\t$ is supported in a much smaller region
$\Omega_\t\subset\mathcal{I}$.
\begin{lem}\label{key_decomp}
If $\text{supp}\widehat{f}\subset\mathcal{I}$ and $1\ll K\ll R$,
then for any $\e>0$ we have the following estimate on the cylinder
$B(0,R)\times [0,R]\subset\R^{2+1}$
\begin{align}
\label{K1} |Tf(x)|&\lesssim R^\e
\max_{\frac{1}{K}\geq\d>\frac{1}{\sqrt{R}}}\max_{\mathcal{E}_\d}
\biggl[\sum_{\Omega_\t\in\mathcal{E}_\d}\Bigl(\psi_\t\,\prod^3_{j=1}|Tf_{\t_j}|^{\frac{1}{3}}\Bigr)^2\biggr]^{1/2}
\\
\label{K2} &\quad\q+R^\e\max_{\mathcal{E}_{\frac{1}{\sqrt{R}}}}
\bigg[\sum_{\Omega_\t\in\mathcal{E}_{\frac{1}{\sqrt{R}}}}\Bigl(\psi_\t\,|Tf_\t|\Bigr)^2\bigg]^{1/2},
\end{align}
where $\^{f}_{\t_j}$ is supported in $\Omega_{\t_j}$ for
$j=1,\,2,\,3$, and
\begin{itemize}
\item $\mathcal{E}_\d$ consists of at most
$\left(\frac{1}\d\right)^{1+\e}$ disjoint $\d\times\d$ squares $\Omega_\t$;
\end{itemize}
\begin{itemize}
\item
$\{\Omega_{\t_j}\}^3_{j=1}$ is
 a triple of non-collinear $\frac{\d}{K}\times \frac{\d}{K}$ squares inside  $\Omega_\t$;
\end{itemize}
\begin{itemize}
\item  For each $\t$, $\psi_\t$ is a non-negative function on $\R^3$ which is constant on unit cubes and satisfies
\begin{equation}\label{aver}
    \frac{1}{|B|}\int_B\psi^4_\t(x)dx\lesssim R^\e,
\end{equation}
for all $B$ taken in a tiling of $\R^3$ with $\mathcal{C}^*_\t-$boxes.
\end{itemize}
\end{lem}

\begin{rem}
A triplet
$\bigl(\Omega_{\t_1},\,\Omega_{\t_2},\,\Omega_{\t_3}\bigr)$ in $\Omega_\t$ is
said to be non-collinear if
\begin{equation}\label{non-cop}
|\xi_{\t_1}-\xi_{\t_2}| \geq|\xi_{\t_1}-\xi_{\t_3}|
\geq\text{dist}\left(\xi_{\t_3},\,\ell(\xi_{\t_1},\xi_{\t_2})\right)
>10^3\frac{\a 2^\a}{K},
\end{equation}
where\ $\xi_j$ is the center of\ $\Omega_{\t_j}$ for $j=1,\,2,\,3$.
Consequently, the caps
$\mathcal{C}_{\t_1},\,\mathcal{C}_{\t_2},\,\mathcal{C}_{\t_3}$ are
transversal, that is, the exterior normal vectors associated to these three caps
are linearly independent, uniformly with respect to the variables belonging to $\Omega_{\t_j}$ for $j=1,2,3$.
We refer to \cite{ref BCT} for the precise description of the transversality
condition, see also Section 3. This condition is ready for the multi-linear
restriction estimate established in \cite{ref BCT}, as frequently
used in \cite{ref Bourgain3,ref Bourgain Guth}.
\end{rem}

\begin{rem}
This lemma is established in the spirit of Bourgain and Guth,
where it differs from \cite{ref Bourgain Guth} in two aspect.
First, the non-collinear condition is reformulated
in \eqref{non-cop} to handle general $\a>1$. Second, the scales of
the caps and dual caps have depended on $\a$ already.
The absence of the algebraic structure of the symbol $|\xi|^\a$ for general
$\a>1$ will lead to some difficulties in the deduction of \eqref{K1} \eqref{K2}
as well as the application of this inequality to the proof of Theorem \ref{thm 1}.
These obstacles make our argument more complicated than  \cite{ref Bourgain3}.
\end{rem}

\begin{rem}\label{local-constant}
If $\^f_{\t}$ is supported in a square $\Omega_\t$ of size $\d$, $|Tf_{\t}|(x)$ can be regarded essentially as
a constant on each $\C^*_\t-$box. We call this the local constancy property
indicated in \cite{ref Bourgain Guth}.
\end{rem}

The advantage of this inequality allows us to gain the
transversality condition in each term of the summation \eqref{K1} by losing
only an $R^\e$ factor. This is favorable especially in proving some
non-endpoint estimates. We also point out that the precise cardinality of $\mathcal{E}_\d$
will not be used in the proof of Theorem \ref{thm 1}. We will prove
Lemma \ref{key_decomp} in Section 5.

\subsection{A primary reduction of the problem}
\noindent
\vskip 0.2cm

By Littlewood-Paley's theory, Sobolev embedding and H\"{o}lder's
inequality, Theorem \ref{thm 1} amounts to showing for any $\e>0$, there is a
$C_{\a,\e}$ such that
\begin{equation}\label{1'}
    \Bigl\|\sup_{0<x_3<1}|Tf(\cdot,x_3)|\Bigr\|_{L^2(B(0,1))}\leq
    C_{\a,\,\e}\, R^{\frac{3}8+\e}\|f\|_2,
\end{equation}
for $ \widehat{f}(\xi)$ supported in $\{\xi\in\R^2\mid R/2\leq
|\xi|\leq2R\}$ with $R$ large enough.

After a re-scaling, we reduce \eqref{1'} to
\begin{equation}\label{2}
\Bigl\|\sup_{0<x_3<R^\a}|Tf(\cdot,x_3)|\Bigr\|_{L^2(B(0,R))} \leq
C_{\a,\e}\,R^{\frac{3}8+\e}\|f\|_2,\q\text{supp}\widehat{f}
\subset\mathcal{I}.
\end{equation}
When $\a=2$, it is observed by Lee \cite{ref Lee} that to get \eqref{2},
it suffices to prove it with the supremum taken only for $0<x_3<R$.
This fact simplifies the problem significantly
so that the result of $s>\frac38$ can be deduced for $d=2$.
This reduction is also necessary for the argument in
\cite{ref Bourgain3}. We extend this result to \emph{all} $\a>1$ by proving
the following lemma.
\begin{lem}\label{lem 2}
Suppose for any $\e>0$, there exists some $C_\e>0$ such that
\begin{equation}\label{lem 2'}
       \Bigl\|
       \sup_{0<x_3<R}|Tf(\cdot,x_3)|\Bigr\|_{L^2(B(0,R))}
       \leq
       C_\e R^{\frac{3}8+\e}\|f\|_2,
\end{equation}
for $R$ sufficient large and $\text{supp}\widehat{f}\subset\{\xi\in
\R^2\mid\frac38\leq|\xi|\leq\frac{17}8 \}$.
Then \eqref{2} holds.
\end{lem}

\begin{rem}
Intuitively, one might expect that the interval over which the
supremum is taken for $x_3$ in \eqref{lem 2'} should be $(0,\,R^\frac{\a}{2})$.
Although this can be deduced easily by modifying our argument slightly,
we will lose more derivatives in Theorem \ref{thm 1} if we use \eqref{lem 2'} with $0<x_3<R^{\frac\a2}$.
The loss of derivatives forces the $s$ in \eqref{1} relies heavily on $\a$ and
this will confine $\a$ in a small range in order to improve Sj\"{o}lin's result.
However, our result can be strengthened so that $s$ is independent of $\a$ thanks to
Lemma \ref{lem 2}. We point out that the global maximal inequality is
$\a-$dependent. See \cite{Rogers08} 
for details.
\end{rem}

\begin{rem}
Heuristically speaking, the idea behind this lemma can be
interpreted in terms of the propagation speed.
If the frequency of the initial data $f$ is
localized at $\mathcal{I}$, then the propagation speed of
$U(t)f$ can be morally regarded as finite. Suppose
$R$ is large enough so that $f$ is mainly concentrated in the ball
$B(0,\,R)$. If one waits at a position in $B(0,\,R)$ for the
maximal amplitude of the solution during the time period $0<t<R^\a$ to occur,
then by the finite speed of propagation, this maximal amplitude can
be expected to happen before the time at $R$.
This heuristic intuition is justified by the following lemma.
\end{rem}

\begin{lem}\label{lem 1}
Let  $\displaystyle \text{supp}\,
\widehat{f}(\xi)\subset\mathcal{I}$ and $j=0,1,\ldots, [R^{\a-1}]$,
$t_j=jR$. Set $I_j=[t_j,\,t_{j+1})$ for $j<[R^{\a-1}]$ and
$I_{[R^{\a-1}]}=[t_{[R^{\a-1}]}, R^\a)$. 
Denote $x=(x',x_{d+1})$ where $x'=(x_1,\ldots,x_{d})$
and take a smooth function
$\varphi\in C^\infty_0(B(0,\,2R))$ so that $\varphi(x')=1$ on the ball
$B(0,\,R)$. Then for any $\e>0$, there is a $C_{\a,\,\e}>0$ and a
family of functions $\{f_j\}_j$ satisfying
$$
\text{supp}\widehat{f}_j\subset
\Bigl\{\xi\in\R^d\,\Big|\;\frac{1}{2}-\frac{1}{R}\leq|\xi|\leq
2+\frac{1}{R}\Bigr\}\mathrel{\mathop=^{\rm
def}}\mathcal{I}_{\frac{1}R}, $$ so that for $x_{d+1}\in I_j$,
\begin{equation}\label{lem 1 '}
       \varphi(x')Tf(x)
       =\varphi(x')\chi_{I_j}(x_{d+1})Tf_j(x',\,x_{d+1}-t_j)
       +\mathcal{O}_{\a,\e}(R^{-99d}\|f\|_2),
\end{equation}
or equivalently, viewing $x_{d+1}=t\in I_j$,
\begin{equation}\label{lem 1 '}
       \varphi(x')U(t)f(x')
       =\varphi(x')\chi_{I_j}(t)U(t-t_j)f_j(x')
       +\mathcal{O}_{\a,\e}(R^{-99d}\|f\|_2),
\end{equation}
Moreover, there exists a positive constant $c_d>0$ such that
\begin{equation}\label{lem 1''}
       \Bigl\|\Bigl(\sum^{[R^{\a-1}]}_{j=0}|f_j|^2\Bigr)^{\frac12}\Bigr\|_2\leq\, C_{\a,\,\e}\, R^{\e c_d} \|f\|_2.
\end{equation}
\end{lem}
To prove this lemma, we introduce a localization argument
which allows us to regard a Schwartz function with compact Fourier frequencies
as a smooth cut-off function by losing $R^\e$. That is  why we have to
lose $R^{\e c_d}$ in \eqref{lem 1''}, but this is adequate for our purpose.

\begin{rem}
In our proof, $f_j$ is constructed by localizing $Tf(x,t_j)$ with Schwartz functions. This leads to a loss of
 $\frac1R-$enlargement of $\mathcal{I}$ in the frequency space, but this does not affect the use of this lemma.
\end{rem}
We end up this section by showing Lemma \ref{lem 2} follows
from Lemma \ref{lem 1}.

\begin{proof}[Proof of Lemma \ref{lem 2}]
 In view of \eqref{lem 1 '}, we have for $d=2$
 $$\varphi(x')|Tf(x',x_3)|\lesssim_{\a,\e}\varphi(x')\sum^{[R^{\a-1}]}_{j=0}|\chi_{I_j}(x_3)Tf_j(x',x_3-t_j)|+R^{-198}\|f\|_2.$$
 Choosing $R$ large enough and neglecting $R^{-198}\|f\|_2$, we obtain
\begin{align*}
 \sup_{0<x_3<R^\a}|\varphi(x')Tf(x)|^2
 \lesssim_{\a,\e} \sum^{[R^{\a-1}]}_{j=0}
 \sup_{0<x_3-t_j<R}|\varphi(x')Tf_j(x',x_3-t_j)|^2.
\end{align*}
 Integrating both sides of the above inequality on $B(0,\,R)$, we may estimate the left side of \eqref{2} by
\begin{align*}
 \Bigl(\sum^{[R^{\a-1}]}_{j=0}\Bigl\|\sup_{0<x_3-t_j<R}|\varphi(x')Tf_j(x',x_3-t_j)|\Bigr\|^2_2\Bigr)^{\frac12}.
\end{align*}
Using \eqref{lem 2'} and \eqref{lem 1''}, we obtain
 $$
 \eqref{2}\lesssim_{\a,\e}R^{\frac{3}{8}+\e}\|f\|_2.
 $$
\end{proof}

\begin{rem}
In proving \eqref{lem 2'}, we always fix an $\e>0$ first and then take $R$
large, which may depend possibly on $\e$, $\a$ and $\|f\|_2$. This allows us to eliminate
as many error terms as possible by repeatedly using the localization argument.
\end{rem}

\begin{rem}
  In fact, the original proof of [L, Lemma 2.1] for the classical Schr\"odinger equation
  works for the generalized case as well. Moreover, the proof of [LR, Lemma 2.1] can be used to get rid of the $\e$ loss.
\end{rem}

\section{Proof of the main result}\label{sect:pf_main}
\noindent

Now we are in the position to prove Theorem \ref{thm 1}. For any fixed $\e>0$,
we normalize $\|f\|_2=1$. In light of Lemma \ref{key_decomp} and \ref{lem 2},
\eqref{lem 2'} amounts
to obtaining the following two estimates for $R$ large enough
\begin{align}
\label{reduced1}
    \sum_{\substack{\frac{1}{\sqrt{R}}<\d<\frac1K\\
    \d\,\text{dyadic}}}&\Bigl[\sum_{\Omega_{\t}:\d\times\d}
    \bigl\|\psi_\t\prod^3_{j=1}|Tf_{\t_j}|^\frac{1}{3}\bigr\|^2_{L^2(|x'|<R)L^{\infty}(|x_3|<R)}\Bigr]^{\frac{1}{2}}
    \lesssim R^{\frac38+\e} ,\\
\label{reduced2}
    \Bigl[&\sum_{\Omega_\t:\frac{1}{\sqrt{R}}\times\frac{1}{\sqrt{R}}}
    \bigl\|\psi_{\t}|Tf_\t|\bigr\|^2_{L^2(|x'|<R)L^{\infty}(|x_3|<R)}\Bigr]^{\frac{1}{2}}
    \lesssim R^{\frac38+\e} ,
\end{align}
where
$x'=(x_1,\,x_2)$
and
$\Omega_\t:\d\times\d$ refers to the partition of $\mathcal{I}_{\frac18}$ into the union of
$\d\times\d-$squares.\\

By orthogonality, it suffices to prove
\begin{equation}\label{reduced3}
    \int_{|x'|<R}\sup_{|x_3|<R}\Bigl|\psi_\t\prod^3_{j=1}|Tf_{\t_j}|^\frac{1}{3}\Bigr|^2(x',x_3)dx'
    \lesssim_\e R^{\frac34+2\e}\|f_\t\|^2_{2},
\end{equation}
and
\begin{equation}\label{reduced4}
    \int_{|x'|<R}\sup_{|x_3|<R}\bigl(\psi_\t|Tf_{\t}|\bigr)^2(x',x_3)dx'\lesssim_\e R^{\frac34+2\e}\|f_\t\|^2_{2},
\end{equation}
where $f_\t$ is defined as $\^f_\t=\^f \chi_{\Omega_\t}$.

\subsection{The proof of \eqref{reduced3}}
\noindent

For brevity, we denote
$G_{\t_1,\t_2,\t_3}=\prod^3_{j=1}|Tf_{\t_j}|^{\frac{1}{3}}$ in
\eqref{reduced3}, where $(\t_1,\,\t_2,\,\t_3)$ corresponds to the squares
$(\Omega_{\t_1},\Omega_{\t_2},\Omega_{\t_3})$
with the properties in Lemma \ref{key_decomp}.
Let $\xi^\t$ be the center of
$\Omega_\t$. Then we may assume  $\xi^\t=(0,\,|\xi^\t|)$ due to the
invariance of \eqref{reduced3} under orthogonal transformations. Let
$\C_\t$ be the $\d-$cap associated to $\Omega_\t$
and $\C^*_\t$ be the polar set of $\C_\t$.
After tiling $B(0,R)\times[0,R]\subset\R^3$ with
$\C^*_\t-$boxes, we have
 $$B(0,R)\times[0,R]\subset\bigcup_{k,j}B_{j,k}$$
where $B_{j,k}$ is a $\C^*_\t-$box labeled by $j$ and $k$,
with $j$ corresponding to the horizontal translation and $k$ to the vertical
(see Figure.1). Adopting the notations in \cite{ref
Bourgain3}, we denote the projection of each $B_{j,k}$ to the
$(x_1,\,x_2)-$variables by $I_j=\pi_{x'}(B_{j,k})$.
Let $\mathcal{P}_{x'}$ be the plane through the point $(x',0)$
and perpendicular to the $x_2-$direction. We define
$J^{x'}_k=\pi_{x_3}(B_{j,k}\cap \mathcal{P}_{x'})$. Then $|J^{x'}_k|\sim\frac1\d$ for all $x'$ and $k$ .
For $I_j$, it is easy to see that the length of the side in direction of $x_1$ is approximately $\frac{1}\d$
and the side in the $x_2-$direction has length $\frac{1}{\d^\a}$.

\begin{figure}[ht]\label{box}
\begin{center}
$$\ecriture{\includegraphics[width=7cm]{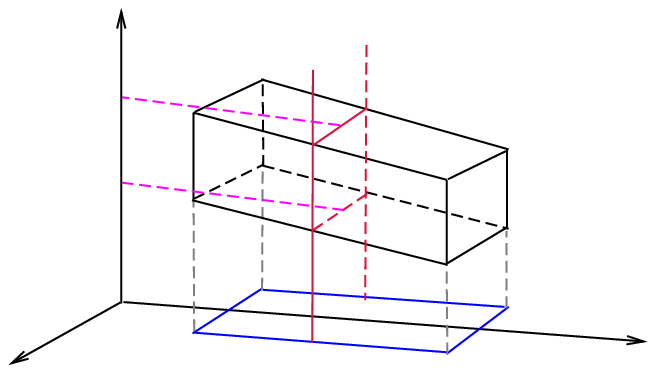}}
{\aat{-1}{-1}{$x_1$}\aat{9}{2}{$O$}\aat{25}{23}{$\mathcal{P}_{x'}$}
\aat{10}{30}{$x_3$}\aat{29}{16}{$B_{j,k}$}\aat{6}{17}{$J^{x'}_k$}\aat{21}{-1}{$I_j=\pi_{x'}(B_{j,k})$}\aat{48}{0}{$x_2$}}$$
\end{center}
\caption{The $\C^*_\t-$box  $B_{j,k}$.}
\end{figure}

By H\"{o}lder's inequality
and the relation $\ell^3\hookrightarrow \ell^\infty$, we have
\begin{align}\label{jia 1}
    &\sum_j\|\psi_\t\G\|^2_{L^2_{x'}(I_j)L^\infty(|x_3|<R)}\nonumber\\
    \leq& \d^{-\frac{1+\a}{3}}\sum_j\|\psi_\t\G\|^2_{L^3_{x'}(I_j)L^\infty(|x_3|<R)}\nonumber\\
    \leq&\d^{-\frac{1+\a}{3}}\sum_j\max_k\|\psi_\t\G\|^2_{L^3_{x'}(I_j)L^\infty_{x_3}(J^{x'}_k)}\nonumber\\
    \leq&\d^{-\frac{1+\a}{3}}\sum_j\Bigl[\sum_k\|\psi_\t\G\|^3_{L^3_{x'}L^\infty_{x_3}(B_{j,k})}\Bigr]^{\frac{2}{3}}.
\end{align}

Using the property \eqref{aver} and Remark \ref{local-constant}, we
have
\begin{align*}
    &\q\|\psi_\t\G\|_{L^3_{x'}L^\infty_{x_3}(B_{j,k})}\\
    &\lesssim\|\psi_\t\|_{L^3_{x'}L^\infty_{x_3}(B_{j,k})}
    \Bigl(\frac{1}{|I_j|}\int_{I_j}\bigl(\frac{1}{|J^{x'}_k|}\int_{J^{x'}_k}|\G(x',x_3)|^3dx_3\bigr)^{\frac{2}{3}}dx'\Bigr)^{\frac{1}{2}}\\
    &\lesssim
    |I_j|^{-\frac{1}{2}}|J_k|^{-\frac{1}{3}}\|\psi_\t\|_{L^3_{x'}L^4_{x_3}(B_{j,k})}\|\G\|_{L^2_{x'}L^3_{x_3}(B_{j,k})}\\
    &\lesssim
    |I_j|^{\frac{1}{3}-\frac{1}{4}-\frac{1}{2}}|J_k|^{-\frac{1}{3}}\|\psi_\t\|_{L^4(B_{j,k})}\|\G\|_{L^2_{x'}L^3_{x_3}(B_{j,k})}\\
    &\lesssim
    \d^{-(1+\a)(\frac{1}{3}-\frac{3}{4})+\frac{1}{3}}R^\e|B_{j,k}|^{\frac{1}{4}}\|\G\|_{L^2_{x'}L^3_{x_3}(B_{j,k})}\\
    &\lesssim R^\e
    \d^{\frac{1}{4}+\frac{\a}{6}}\|\G\|_{L^2_{x'}L^3_{x_3}(B_{j,k})},
\end{align*}
where, in the second inequality, we used the fact that $\psi_\t$ is constant on the unit cubes.
This allows us to control the $L^\infty$ norm of $\psi_\t$ with respect to the $x_3$ variable on $J^{x'}_k$
by $L^4$ norm.
\vskip0.2cm

Plugging this into \eqref{jia 1}, we get by Minkowski's inequality
\begin{align}
   \nonumber |\eqref{jia 1}|&\lesssim \d^{-\frac{1+\a}{3}}\sum_j\Bigl[\sum_k
    R^{3\e}\d^{\frac{3}{4}+\frac{\a}{2}}\|\G\|^3_{L^2_{x'}L^3_{x_3}(B_{j,k})}\Bigr]^{\frac{2}{3}}\\
   \nonumber &\lesssim
    R^{2\e}\d^{\frac{1}{6}}\sum_j\int_{I_j}\Bigl(\sum_k\int_{J^{x'}_k}|\G|^3(x',x_3)dx_3\Bigr)^{\frac{2}{3}}dx'\\
   \label{crux} &\lesssim
    R^{2\e}\d^{\frac{1}{6}}\|\G\|^2_{L^2(|x'|<R)L^3(|x_3|<R)}.
\end{align}
\begin{rem}
In \eqref{crux}, $\a$ is absent from the exponent of
$\d$ and $R$. However, the role of $\a$ is enrolled in the estimation of
\beq\label{G}\|\G\|^2_{L^2(|x'|<R)L^3(|x_3|<R)}.\eeq To handle this
expression, Bourgain used Galilean's transformation to shift the center of
the domain for the integral $Tf_\nu$ to the origin when $\a=2$. We can not directly use this
due to the absence of the algebraic structure of $|\xi|^\a$ for general $\a>1$.
To adapt his strategy, we get around this obstacle by using Taylor's
expansion.
We also use a localization argument as in the proof of Lemma \ref{lem 1}.
\end{rem}

It remains  to evaluate \eqref{G}. Let us introduce some notations. Define
$$
\T_{\d,\,\t}h(x_1,\,x_2,\,x_3)=\int_{\R^2}e^{i[x_1\eta_1+x_2\eta_2+x_3\Phi\,(\xi^\t,\ld_0,\a,\d,\eta)]}\chi(\eta)\^{h}(\eta)d\eta,
$$
where $\chi(\eta)$ is a smooth function adapted to the unit square $\Omega$ and
$$\Phi(\xi^\t,\ld_0,\a,\d,\eta)=\frac\a2|\xi^\t|^{\a-2}\big(\eta^2_1+(\a-1)\eta^2_2\big)+\Theta(\eta)\d|\eta|^3,$$
\beq\nonumber \Theta(\eta)
=\frac{\a(\a-2)}{6}|\xi^\t+\d\eta\ld_0|^{\a-3}
\Big[3\langle\overline{\xi^\t+\d\eta\ld_0},\overline{\eta}\rangle+(\a-4)\langle\overline{\xi^\t+\d\eta\ld_0},\overline{\eta}\rangle^3\Big],
\eeq with $\ld_0\in (0,1)$.\\

Let
$\widehat{g^\d_{\t_\nu}}(\eta)=\d\widehat{f_{\t_\nu}}(\xi^\t+\d\eta)$
and consider Taylor's expansion of $|\xi^\t+\d\eta|^\a$ at $\xi^\t$
up to the third order. We have for some
$\ld_0\in(0,1)$
$$
|Tf_{\t_\nu}(x_1,\,x_2,\,x_3)|=\d\Big|\T_{\d,\,\t}(g^\d_{\t_{\nu}})\Bigl(\d x_1,\,\d(
x_2+x_3\a|\xi^\t|^{\a-1}),\,\d^2 x_3\Bigr)\Big|.
$$
Using H\"{o}lder's inequality in $x_2$ then making change of
variables as
$$(x_1,\,x_2,\,x_3)\rta(\d^{-1}y_1,\,\d^{-1}y_2,\,\d^{-2}y_3),$$
we get
\begin{align*}
&\|\G\|^2_{L^2(|x'|<R)L^3(|x_3|<R)}
\\
\lesssim&\int\limits_{|x'|< R}\Big(\int\limits_{|x_3|<
R}\prod^3_{\nu=1}\d\Big|\T_{\d,\t}(g^\d_{\t_\nu})\big(\d x_1,\d(x_2+\a
x_3|\xi^\t|^{\a-1}),\d^2
x_3\big)\Big|dx_3\Big)^{\frac{2}{3}}dx'
\\
\lesssim&\d^2 R^{\frac{1}{3}}\int\limits_{|x_1|<2R}\Bigg(\iint\limits_{\substack{|x_3|<R\\|x_2|<\a
2^\a
R}}\prod^3_{\nu=1}\Big|\T_{\d,\t}(g^\d_{\t_\nu})(\d
x_1,\d x_2,\d^2 x_3)\Big|dx_2dx_3\Bigg )^{\frac{2}{3}}dx_1
\\
\lesssim&\d^2 R^{\frac{1}{3}}\int\limits_{|y_1|<2\d R} \Bigg(\iint\limits_{\substack{|y_3|<\d^2
R\\|y_2|<\a
2^\a \d R}}\prod^3_{\nu=1}\Big|\T_{\d,\t}(g^\d_{\t_\nu})(y_1,y_2,y_3)\Big|\frac{dy_2}{\d}\frac{dy_3}{\d^2}\Bigg)^{\frac{2}{3}}\frac{dy_1}{\d}.
\end{align*}
Partitioning the range for $y_1$  into consecutive intervals $I_\mu$ as follows
$$
(-2\d R,\,2\d R)=\bigcup_\mu I_\mu\,,\q|I_\mu|=\d^2 R,
$$
we get
\begin{align*}
&\|\G\|^2_{L^2(|x'|<R)L^3(|x_3|<R)}\\
\lesssim&R^{\frac{1}{3}}\d^{-1}\sum_{\mu}\int_{I_\mu}
\Big\|\prod^3_{\nu=1}\T_{\d,\t}(g^\d_{\t_\nu})(y_1,\cdot,\cdot)\Big\|^{\frac{2}{3}}_{L^1(|y_2|<\a
2^\a \d R; |y_3|<\d^2 R)}dy_1.
\end{align*}
Applying H\"{o}lder's inequality with respect to $y_1$ on each $I_\mu$ and then decomposing the interval for $y_2$ similarly as
$$
(-2\d R,\,2\d R)=\bigcup_{\mu'} I_{\mu'}\,,\q\;|I_{\mu'}|=\d^2 R,
$$
we obtain
\begin{align}
\nonumber&\|\G\|^2_{L^2(|x'|<R)L^3(|x_3|<R)}\\
\nonumber\lesssim&
R^{\frac{2}{3}}\d^{-\frac{1}{3}}\sum_\mu\Big(\sum_{\mu'}
\int\limits_{y_2\in
I_{\mu'}}dy_2\iint\limits_{\substack{y_1\in I_\mu\\|y_3|<\d^2
R}}\prod^3_{\nu=1}\Big|\T_{\d,\t}(g^\d_{\t_\nu})(y_1,y_2,y_3)\Big|
dy_1 dy_3\Big)^{\frac{2}{3}}
\\
\label{3.6'}\lesssim&R^{\frac{2}{3}}\d^{-\frac{1}{3}}\sum_{\mu,\,\mu'}\Big\|\prod^3_{\nu=1}\T_{\d,\t}(g^\d_{\t_\nu})\Big\|^{\frac{2}{3}}_{L^1(Q_{\mu,\mu'}\times[-\d^2
R,\d^2 R])},
\end{align}
where $Q_{\mu,\mu'}=I_\mu\times I_{\mu'}$ is a $\d^2 R\times \d^2 R-$square.
\vskip0.2cm
 To evaluate
$$
\Bigl\|\prod^3_{\nu=1}\T_{\d,\t}(g^\d_{\t_\nu})\Bigr\|_{L^1(Q_{\mu,\mu'}\times[\d^2R,\,\d^2 R])},
$$
we need to introduce a localization argument based on Poisson summation, with respect to the $(y_1,\,y_2)-$variables.
\vskip0.2cm
Denote the center of $Q_{\mu,\,\mu'}$ by $y'_{\mu,\mu'}$, which belongs to $\displaystyle \d^{2} R\,\Z^2\mathrel{\mathop=^{\rm def}}\mathcal{Z}.$
Choose a Schwartz function $\beta\geq0$ such that
supp $\widehat{\beta}\subset B(0,1/2)\subset\R^2$ and $\^\beta(0)=1$. We have for all $z\in\R^2$
\begin{equation}\label{betapoisson}
\d^{2\e}\sum_{y'_{\mu,\mu'}\in\mathcal{Z}}\beta\Bigg(\frac{y'_{\mu,\mu'}-z}{\d^{2-\e}R}\Bigg)=1.
\end{equation}
Fix $y'_{\mu_0,\,\mu'_0}\in\mathcal{Z}$ and define
$$K(y',\,y_3,\,z)=\int_{\R^2_\eta} e^{i[\langle y'-z,\,\eta\rangle+y_3\Phi(\xi^\t,\,\ld_0,\,\a,\,\d,\,\eta)]}\chi_{Q_{\mu_0,\mu'_0}}(y')\chi(\eta)d\eta.$$
Then
\begin{align}\label{localozation split}
\Big|\chi_{Q_{\mu_0,\mu'_0}}(y')\T_{\d,\t}(g^\d_{\t_\nu})(y',y_3)\Big|
\lesssim F_1+F_2,
\end{align}
where
\begin{align*}
F_1=&\d^{2\e}\Big|\int K(y',y_3,z)
\sum_{\substack{y'_{\mu,\mu'}\in\mathcal{Z}\\|y'_{\mu,\mu'}-y'_{\mu_0,\mu'_0}|\leq\d^{2}R^{1+\e}}}
\beta\Bigg(\frac{y'_{\mu,\mu'}-z}{\d^{2-\e}R}\Bigg)g^\d_{\t_{\nu}}(z)dz\Big|,\\
F_2=&\d^{2\e}\Big|\int K(y',y_3,z)
\sum_{\substack{y'_{\mu,\mu'}\in\mathcal{Z}\\|y'_{\mu,\mu'}-y'_{\mu_0,\mu'_0}|>\d^{2}R^{1+\e}}}
\beta\Bigg(\frac{y'_{\mu,\mu'}-z}{\d^{2-\e}R}\Bigg)g^\d_{\t_{\nu}}(z)dz\Big|.\\
\end{align*}

First, we estimate $F_2$  in the following manner
$$F_2\leq\mathfrak{F}_{2,1}+\mathfrak{F}_{2,2},$$
where
\begin{align*}
\mathfrak{F}_{2,1}=\d^{2\e}\int\limits_{\substack{|z-y'_{\mu_0,\mu'_0}|\\ \leq \a 2^{4\a}\d^2R^{1+\e_1}}} |K(y,z)|
\sum_{\substack{z_{\mu,\mu'}\in\Z^2\\ \Big|z_{\mu,\mu'}-\frac{y'_{\mu_0,\mu'_0}}{\d^{2}R}\Big|>R^{\e}}}
\beta\Big(\d^\e\bigl(z_{\mu,\mu'}-\frac{z}{\d^{2}R}\bigr)\Big)|g^\d_{\t_{\nu}}|(z)dz,\\
\mathfrak{F}_{2,2}=\d^{2\e}\int\limits_{\substack{|z-y'_{\mu_0,\mu'_0}|\\ > \a 2^{4\a}\d^2R^{1+\e_1}}} |K(y,z)|
\sum_{\substack{z_{\mu,\mu'}\in\Z^2\\ \Big|z_{\mu,\mu'}-\frac{y'_{\mu_0,\mu'_0}}{\d^{2}R}\Big|>R^{\e}}}
\beta\Big(\d^\e\bigl(z_{\mu,\mu'}-\frac{z}{\d^{2}R}\bigr)\Big)|g^\d_{\t_{\nu}}|(z)dz,
\end{align*}
with $z_{\mu,\mu'}=y'_{\mu,\mu'}\d^{-2}R^{-1}$ and $\e_1=0.01\e$.
\vskip0.2cm

Since $R$ can be chosen large enough so that $R^\e\gg \a 2^{4\a}>1$, we have in $\mathfrak{F}_{2,1}$
$$
\bigl|z_{\mu,\mu'}-\frac{z}{\d^2R}\bigr|
\geq\Bigl|z_{\mu,\mu'}-\frac{y'_{\mu_0,\mu'_0}}{\d^{2}R}\Bigr|-\Bigl|\frac{z}{\d^2R}-\frac{y'_{\mu_0,\mu'_0}}{\d^{2}R}\Bigr|
\geq\frac{R^{0.9\e}}{2}.
$$
Hence
\begin{equation}\label{outbetasum}
\sum_{ \Big|z_{\mu,\mu'}-\frac{y'_{\mu_0,\mu'_0}}{\d^{2}R}\Big|>R^{\e}}
\d^{2\e}\beta\Big(\d^\e\bigl(z_{\mu,\mu'}-\frac{z}{\d^{2}R}\bigr)\Big)
\end{equation}
is bounded by
\begin{align*}
\d^{2\e}\int_{|z|>\frac{R^{0.9\e}}{2}}\beta(\d^\e z)dz
\lesssim_N &\int_{|z|>0.5(\d R^{0.9})^\e}(1+|z|)^{-N}dz.
\end{align*}
Noting that $\d>R^{-\frac12}$, we have for suitably large $N$ depending on $\e$
$$
\eqref{outbetasum}\lesssim_\e R^{-2000}.
$$
By Cauchy-Schwarz's inequality in $z-$variables and the boundedness of $\|K(y,\cdot)\|_2$, we obtain
$$
\mathfrak{F}_{2,1}\lesssim_\e R^{-2000}\int \bigl|K(y',y_3,z)g^\d_{\t_{\nu}}(z)\bigr|dz\lesssim_\e R^{-2000}\|g^\d_{\t_\nu}\|_2.
$$

To estimate $\mathfrak{F}_{2,2}$, we write in view of \eqref{betapoisson}
\begin{equation}\label{estF22}
\mathfrak{F}_{2,2}\leq\int_{|z-y'_{\mu_0,\mu'_0}| > \a 2^{4\a}\d^2R^{1+\e_1}} |K(y',y_3,z)|\cdot
|g^\d_{\t_{\nu}}|(z)dz.
\end{equation}
Since $y'$ is restricted in a $2\d^2R-$nieghborhood of $y'_{\mu_0,\mu_0'}$,
we have
\begin{align*}
 &|y'-z|-|y_3\nabla_{\eta}\Phi(\xi^\t,\ld_0,\a,\d\eta)|\\
\geq&|z-y'_{\mu_0,\mu'_0}|-|y'-y'_{\mu_0,\mu'_0}|-|y_3\nabla_{\eta}\Phi(\xi^\t,\ld_0,\a,\d\eta)|\\
\gtrsim&\a2^{4\a}\d^2 R^{1+\e_1}-2\d^2R-\a 2^{3\a}\d^2R\\
\gtrsim&\a 2^{\a-1}\d^2R^{1+\e_1}.
\end{align*}
By introducing the differential operator
$$
\mathfrak{D}=\frac{y'-z+y_3\nabla_\eta\Phi}{|y'-z+y_3\nabla_\eta\Phi|^2}\cdot\nabla_{\eta},
$$
we may estimate $K(y',y_3,z)$ in $\mathfrak{F}_{2,2}$ using integration by parts to get
$$
K(y',y_3,z)\lesssim_N |y'-z|^{-N}.
$$
Inserting this to \eqref{estF22} and using Cauchy-Schwarz, we have for suitable $N=N(\e_1)$
$$\mathfrak{F}_{2,2}\lesssim_{\a,\e}R^{-2000}\|g^\d_{\t_{\nu}}\|_2.$$
Consequently, the contribution of $F_2$ to \eqref{localozation split} is negligible.
\vskip 0.2cm

Now, let us evaluate the contribution of $F_1$. For brevity, we denote
$$B_{\mu_0,\mu'_0}(z)=\sum_{\substack{y'_{\mu,\mu'}\in\mathcal{Z}\\
|y'_{\mu,\mu'}-y'_{\mu_0,\mu'_0}|\leq\d^{2}R^{1+\e}}}\beta\Bigg(\frac{y'_{\mu,\mu'}-z}{\d^{2-\e}R}\Bigg).$$
From the definition of $\widehat{g^\d_{\t_\nu}}$, we have
$$
\text{supp
}\widehat{B_{\mu_0,\mu'_0}g^\d_{\t_\nu}}\subset\frac{1}{\d}(\Omega_{\t_\nu}-\xi^\t)+\mathcal{O}(\frac{1}{R^{\frac\e2}}).
$$
Note that we may choose $K\ll R^{\frac\e2}$, the non-collinear condition is still fulfilled by
the supports of
$\{\widehat{B_{\mu_0,\mu'_0}g^\d_{\t_\nu}}\}^3_{\nu=1}$. In the
meantime, the main contribution of
$$
\|\prod^3_{\nu=1}\T_{\d,\t}(g^\d_{\t_\nu})\|_{L^1(Q_{\mu,\,\mu'}\times[-\d^2
R,\,\d^2 R])}
$$
comes from
\begin{align}\label{hhh}
\d^{2\e}\Big\|\prod^3_{\nu=1}\T_{\d,\t}(B_{\mu,\,\mu'}g^\d_{\t_\nu})\Big\|_{L^1(Q_{\mu,\mu'}\times[-\d^2
R,\,\d^2 R] )}.
\end{align}
For the $\e>0$ at the beginning of this section, we claim
\begin{equation}\label{claimbct}
\Big\|\prod^3_{\nu=1}\T_{\d,\t}(B_{\mu,\mu'}g^\d_{\t_\nu})\Big\|_{L^1(Q_{\mu,\mu'}\times[-\d^2
R,\,\d^2 R])}\lesssim_\e
R^\e\prod^3_{\nu=1}\|B_{\mu,\mu'}g^\d_{\t_\nu}\|_{L^2}.
\end{equation}
We postpone the proof of \eqref{claimbct} to Subsection 3.3. At
present, we show how \eqref{claimbct} implies \eqref{reduced3}.
Using \eqref{3.6'}, \eqref{claimbct} and H\"{o}lder's inequality, we obtain
$$
\eqref{crux}\lesssim
R^{\frac{2}{3}+3\e}\d^{\frac{1}{6}-\frac{1}{3}}
\prod^3_{\nu=1}\Big(\sum_{\mu,\mu'}\|B_{\mu,\mu'}g^\d_{\t_\nu}\|^2_{L^2}\Big)^{\frac{1}{3}}.
$$
From the definition of $B_{\mu.\mu'}$, we have by
Cauchy-Schwarz
\begin{align}
\label{local1}&\sum_{\mu,\mu'}\int |B_{\mu,\mu'}(z)g^\d_{\t_\nu}(z)|^2dz\\
\nonumber\leq&
R^{2\e}\sum_{\mu,\mu'}\sum_{\substack{y'_0\in\mathcal{Z}\\|y'_0-y'_{\mu,\mu'}|\leq
\d^2 R^{1+\e}}}
\int_{\R^2}\Big|\beta\Big(\frac{y'_0-z}{\d^{2-\e}R}\Big)g^\d_{\t_\nu}(z)\Big|^2dz\\
\nonumber\leq&R^{2\e}\sum_{y'_0\in\mathcal{Z}}r_{y'_0}\int_{\R^2}\Big|\beta\Big(\frac{y'_0-z}{\d^{2-\e}R}\Big)g^\d_{\t_\nu}(z)\Big|^2dz,
\end{align}
where
$$r_{y'_0}=\#\bigl\{y'_{\mu,\mu'}\in\mathcal{Z}\mid
\bigl|y'_0-y'_{\mu,\mu'}\bigr|\leq \d^2 R^{1+\e}\bigr\}\lesssim R^{2\e}.$$
Invoking $\frac{1}{\d}<\sqrt{R}$, we get
\begin{align*}
\eqref{local1}\lesssim&
R^{4\e}\d^{-4\e}\int\Big[\d^{2\e}\sum_{y'_0\in\mathcal{Z}}\beta\Big(\frac{y'_0-z}{\d^{2-\e}R}\Big)|g^\d_{\t_\nu}(z)|\Big]^2dz\\
\lesssim&\left(\frac{R}{\d}\right )^{4\e}\|g^\d_{\t_\nu}\|^2_2\\
\lesssim& R^{6\e}\|f_\t\|^2_2.
\end{align*}
As a consequence, we have
$$
\eqref{crux}\lesssim_\e\, R^{\frac34+9\e}\|f_\t\|^2_2.
$$
This implies
\eqref{reduced3} since $\e>0$ can be taken arbitrarily small.

\subsection{The proof of \eqref{reduced4}}
\noindent
\vskip 0.2cm

Letting $\d=\frac{1}{\sqrt{R}}$, we adopt
the same argument as in Subsection 3.1 to obtain
\eqref{crux} with $T f_\t$ in place of $\G$ such that
$$\int_{|x'|<R}\sup_{|x_3|<R}\bigl(\psi_\t|Tf_{\t}|\bigr)^2(x',x_3)dx'\lesssim_\e R^\e\d^{\frac16}\|Tf_\t\|^2_{L^2(|x'|<R)L^3(|x_3|<R)},$$
 where $\Omega_{\t}$ is a
$\frac{1}{\sqrt{R}}\times\frac{1}{\sqrt{R}}-$square.

Denoting
$\widehat{g^\d_{\t}}(\eta)=\d\widehat{f_\t}(\xi^\t+\d\eta)$ and
performing the previous argument, we have
$$
\int_{|x'|<R}\sup_{|x_3|<R}\bigl(\psi_\t|Tf_{\t}|\bigr)^2(x',x_3)dx'\lesssim_\e
R^{\frac{1}{12}+\e}\sum_{\mu,\,\mu'}\|\mathcal{T}_{\d,\t}(g^\d_{\t})\|^2_{L^3(Q_{\mu,\,\mu'}\times [-1,1])},
$$
where $Q_{\mu,\,\mu'}$ is a square with unit length. Invoking the definition of $\mathcal{T}_{\d,\t}
(g^{\d}_{\t})$, we have
$$\|\mathcal{T}_{\d,\t}
(g^{\d}_{\t})\|_{L^\infty(Q_{\mu,\mu'}\times[-1,1])}\leq\|g^{\d}_{\t}\|_2\leq\|f_\t\|_2.$$
By H\"{o}lder's inequality and Plancherel's theorem, we obtain
\begin{align*}
    &\sum_{\mu,\mu'}\|\mathcal{T}_{\d,\t}(g^\d_{\t})\|^2_{L^3(Q_{\mu,\mu'}\times[-1,1])}\\
    \lesssim& \|f_\t\|^{\frac23}_2\sum_{\mu,\mu'}\|\mathcal{T}_{\d,\t}(g^\d_{\t})\|^{\frac43}_{L^2(Q_{\mu,\mu'}\times[-1,1])}\\
    \lesssim&\;\|f_\t\|^{\frac23}_2\bigl(\sum_{\mu,\mu'}1\bigr)^{\frac13}
    \Bigl(\sum_{\mu,\mu'}\|\mathcal{T}_{\d,\t}(g^\d_{\t})\|^2_{L^2(Q_{\mu,\mu'}\times[-1,1])}\Bigr)^{\frac23}\\
    \lesssim&\;R^{\frac23}\|f_{\t}\|^{\frac23}_2\|\mathcal{T}_{\d,\t}(g^\d_{\t})\|^{\frac43}_{L^\infty(|x_3|<1)L^2(|x'|\leq R)}\\
    \lesssim&\;R^{\frac23}\|f_\t\|^2_2.
\end{align*}
Therefore \eqref{reduced4} follows. Collecting \eqref{reduced3} and \eqref{reduced4},
we conclude that \eqref{claimbct} implies \eqref{lem 2'}. We shall prove \eqref{claimbct} in the next subsection.

\subsection{The proof of \eqref{claimbct}}
\noindent
\vskip 0.2cm

To prove \eqref{claimbct}, we need the multilinear restriction theorem in \cite{ref BCT}.
Since a special form of this theorem is adequate
for our purpose, we formulate it only in this form
and one should consult \cite{ref BCT} for the general statement.

Now we introduce some base assumptions. Let $U\subset\R^{d-1}_\eta$ be a compact neighborhood of the origin
and $\Sigma:U\rightarrow \R^d$ be a smooth parametrisation of a
$(d-1)$ hypersurface of $\R^d$. For $U_\nu\subset U$
and $g_\nu$ supported in $U_\nu\subset\R^{d-1}$ with $1\leq \nu\leq d$, assume that
there is a constant $\mu>0$ such that,
\begin{equation}\label{non-co-condition}
\text{det}\Bigl(X(\eta^{(1)}),\ldots,X(\eta^{(d)})\Bigr)>\mu
\end{equation}
for all $\eta^{(1)}\in U_1,\ldots,\eta^{(d)}\in U_d$, where
$$
X(\eta)=\bigwedge^{d-1}_{k=1}\frac{\partial}{\partial
\eta_k}\Sigma(\eta),\; \eta=(\eta_1,\ldots,\eta_{d-1}).
$$
Assume also that there is a constant
$A\geq 0$
\begin{equation}\label{smooth-condi}
\|\Sigma\|_{C^2(U_\nu)}\leq
A\,,\q\text{for all}\;1\leq \nu\leq d.
\end{equation}
For each $1\leq \nu\leq d$, define for $g_\nu\in L^p(U_\nu)$, $p\geq1$
$$
\mathcal{S} g_\nu(x)=\int_{U_\nu}e^{ix\cdot\Sigma(\eta)}g_\nu(\eta)d\eta.
$$

\begin{thm}\label{bct}
Under the assumption of \eqref{non-co-condition}
and \eqref{smooth-condi}, we have for each $\e>0$,
$q\geq\frac{2d}{d-1}$ and $p'\leq \frac{d-1}{d}q$, there is a
constant $C>0$, depending only on $A,\e,p,q,d,\mu$, for which
\begin{equation}\label{BCT}
\Bigl\|\prod^d_{\nu=1}\mathcal{S} g_\nu\Bigr\|_{L^{\frac qd}(B(0,R))}\leq
C_\e R^\e\prod^d_{\nu=1}\|g_\nu\|_{L^p(U_\nu)},
\end{equation}
for all $g_1,\ldots,g_d\in L^p(\R^{d-1})$ and all $R\geq 1$.
\end{thm}

\begin{rem}
 We shall use \eqref{BCT} below with $d=q=3$ and $p=2$.
\end{rem}

The proof of \eqref{claimbct}  amounts to show
\begin{equation}\label{variablebct}
\Big\|\prod^3_{\nu=1}|\T_{\d,\t}(g^\d_{\t_\nu})|\Big\|_{L^1(B(0,\ld))}\leq
C_\e\ld^\e\prod^3_{\nu=1}\|g^\d_{\t_{\nu}}\|_{L^2},\;\forall\,\ld>0.
\end{equation}
To prove \eqref{variablebct}, we use \eqref{BCT} with
$$\mathcal{S}=\T_{\d,\t},$$
$$U_\nu=\d^{-1}(\Omega_{\t_\nu}-\xi^{\t})+\mathcal{O}(\frac{1}{R^{\frac\e2}}),\,\nu=1,2,3$$
$$g_\nu=\^{g^\d_{\t_\nu} }$$  and
$$
\Sigma:(\eta_1,\eta_2)\rightarrow \bigl(\eta_1,\eta_2,\Phi(\xi^\t,\,\ld_0,\,\a,\,\d,\,\eta)\bigr).
$$
If $\Sigma$ satisfies \eqref{non-co-condition} and \eqref{smooth-condi},
then \eqref{variablebct} follows immediately. Since
the smoothness condition \eqref{smooth-condi} is clear from the definition of $\Phi$,
we only need to show the transversality condition \eqref{non-co-condition}.\\

A simple calculation yields
$$
\begin{cases}
 &\partial_{\eta_1}\Sigma=(1,\,0,\,\partial_{\eta_1}\Phi(\xi^\t,\,\ld_0,\,\a,\,\d,\,\eta)\,)\\
 &\partial_{\eta_2}\Sigma=(0,\,1,\,\partial_{\eta_2}\Phi(\xi^\t,\,\ld_0,\,\a,\,\d,\,\eta)\,)
\end{cases}
$$
where
$$
\begin{cases}
 &\partial_{\eta_1}\Phi=\a|\xi^\t|^{\a-2}\eta_1+\d\partial_{\eta_1}\big(\Theta(\eta)|\eta|^3\big)\\
 &\partial_{\eta_2}\Phi=\a(\a-1)|\xi^\t|^{\a-2}\eta_2+\d\partial_{\eta_2}\big(\Theta(\eta)|\eta|^3\big).
 \end{cases}
$$
Since  \beq\nonumber \Theta(\eta)
=\frac{\a(\a-2)}{6}|\xi^\t+\d\eta\ld_0|^{\a-3}
\Big[3\langle\overline{\xi^\t+\d\eta\ld_0},\overline{\eta}\rangle+(\a-4)\langle\overline{\xi^\t+\d\eta\ld_0},\overline{\eta}\rangle^3\Big],
\eeq we have
$$\nabla_\eta\big(\Theta(\eta)|\eta|^3\big)=\mathcal{O}_\a(1).$$
This along with $\d\leq \frac 1K\ll 1$ allows us to write
$$
X(\eta)
=\Bigl(\a(\a-1)|\xi^\t|^{\a-2}\eta_1,\,\a|\xi^\t|^{\a-2}\eta_2,\,-1\Bigr)+\mathcal{O}_\a(1)\Bigl(\frac{1}{K},\,\frac{1}{K},\,0\Bigr),$$
hence we have
\begin{align}
&\nonumber\text{det}\Bigl(X(\eta^{\t_1}),\,X(\eta^{\t_2}),\,X(\eta^{\t_3})\Bigr)\\
\label{trans}=&\a^2(\a-1)|\xi^\t|^{2(\a-2)} \text{det} \left(
\begin{array}{ccc}
-1& -1 & -1  \\
\eta^{\t_1}_1& \eta^{\t_2}_1 & \eta^{\t_3}_1\\
\eta^{\t_1}_2& \eta^{\t_2}_2 & \eta^{\t_3}_2
\end{array}
\right) +\mathcal{O}_\a\Bigl(\frac{1}{K}\Bigr).
\end{align}
In view of the non-collinear condition fulfilled  by
$\Omega_{\t_1},\,\Omega_{\t_2},\,\Omega_{\t_3}$, we see the area of
the triangle formed by $\eta^{\t_1}$, $\eta^{\t_2}$ and
$\eta^{\t_3}$ is uniformly away from zero, or equivalently, there is a $C>0$ such that
$$
\left| \text{det}\Biggl(
\begin{array}{ccc}
  -1& -1& -1\\
\eta^{\t_1}_1& \eta^{\t_2}_1& \eta^{\t_3}_1\\
\eta^{\t_1}_2& \eta^{\t_2}_2& \eta^{\t_3}_2
\end{array}
\Biggr)\right| \geq C
$$
for all
$\eta^{\t_\nu}\in U_\nu,\,\nu=1,2,3$. Therefore,
we can rearrange the order of the columns in \eqref{trans} to ensure
$$
\text{det}\Biggl(
\begin{array}{ccc}
  -1& -1& -1\\
\eta^{\t_1}_1& \eta^{\t_2}_1& \eta^{\t_3}_1\\
\eta^{\t_1}_2& \eta^{\t_2}_2& \eta^{\t_3}_2
\end{array}
\Biggr)\geq C>0.
$$
Next, we take $K$ large enough so that
$$\eqref{trans}\geq\a^2(\a-1)|\xi^\t|^{\a-2}\frac C2>0.$$
Consequently, we have \eqref{variablebct} and this completes the proof of Theorem
\ref{thm 1}.

\section{Proof of Lemma \ref{lem 1}}\label{sec:lem 1}
\noindent

In this section, we prove Lemma \ref{lem 1}.
Take $\eta\in\mathcal{S}(\R^d)$ such that
$\eta\geq0$ and $\widehat{\eta}$ is supported in the ball
$B(0,\,1/2)$ with $\widehat{\eta}(0)=1$.
Denoting by $\mathfrak{X}=R\,\Z^d$, we have Poisson's summation formula
$$
\sum_{x^\sigma\in\mathfrak{X}} \eta\Bigl(\frac{x-x^\sigma}{R}\Bigr)
=1,\q\forall\, x\in\R^d.
$$
We adopt the notions in Lemma \ref{lem 1}.
First, noting that
$$
U(t)f(x')=U(t-t_j)U(t_j)f(x'),
$$
we may write, for $x_{d+1}\in (0,R^\a)$
\begin{equation}\label{2.1}
\psi(x')Tf(x)
=\sum^{[R^{\a-1}]}_{j=0}\int_{\R^{d}}\chi_{I_j}(x_{d+1})\,K(x,y; t_j)\,Tf(y,t_j)\,dy,
\end{equation}
with $K(x,y; t_j)$ defined by
$$
 K(x,y; t_j)=\int_{\R^{d}}
 e^{i[(x'-y)\cdot\xi+(x_{d+1}-t_j)|\xi|^\a]}\,\psi(x')\,\chi(\xi)\,d\xi,
$$
where $\chi\in C^\infty_c(\R^d)$ such that $\chi(\xi)\equiv 1$ for
$\xi\in\mathcal{I}_{\frac1R}$.
Without loss of generality, we may assume $x_{d+1}$ belongs to $ I_j$ for some
$j\in\{0,\ldots,\,[R^{\a-1}]\}$. Using Poisson's summation formula, we have
$$\eqref{2.1}
  =J_1+J_2$$
  where
\begin{align*}
&J_1=\sum_{ \substack{y^\sigma\in\mathfrak{X},\\
  |y^\sg|\geq 10 R^{1+\e}}}\int_{\R^d}K(x,y;t_j)\,\eta\Bigl(\frac{y-y^\sg}{R}\Bigr)\,Tf(y,t_j)dy\\
&J_2=\sum_{ \substack{y^\sigma\in\mathfrak{X},\\
  |y^\sg|< 10 R^{1+\e}}}
  \int_{\R^d}K(x,y;t_j)\,\eta\Bigl(\frac{y-y^\sg}{R}\Bigr)\,Tf(y,t_j)dy\\
\end{align*}

\subsection{ The estimation of $J_1$.}
 Divide the integral $J_1$ into two parts as follows
$$J_1\leq J_{1,1}+J_{1,2},$$
where
\begin{align*}
 J_{1,1}\mathrel{\mathop=^{\rm def}}&\int_{|y|\leq \a 2^{\a+2}R}\bigl|K(x,y;t_j)\bigr|
 \sum_{ y^\sigma\in\mathfrak{X}, |y^\sg|\geq 10 R^{1+\e}} \eta\Bigl(\frac{y-y^\sg}{R}\Bigr)|Tf(y,t_j)|\,dy,\\
 J_{1,2}\mathrel{\mathop=^{\rm def}}&\int_{|y|>\a 2^{\a+2} R}\bigl|K(x,y;t_j)\bigr|\,|Tf(y,t_j)|dy.
\end{align*}
We show $J_1$ is negligible by estimating the contribution of $J_{1,1}$ and $J_{1,2}$,
separately.

\vskip0.4cm

$\bullet$ \textbf{ Estimation of $J_{1,1}$}. It is easy to see that
\begin{align*}
|J_{1,1}|\leq\int_{|y|\leq \a 2^{\a+2}R}
\sum_{z\in\mathbb Z^d,|z|\geq10R^\e}\eta\Bigl(\frac{y}{R}-z\Bigr)\bigl|K(x,y;t_j)\bigr|\,|Tf(y,t_j)|\,dy.
\end{align*}
Since $R^\e\gg \a2^{\a+1}>1$, $|y|\leq \a2^{\a+2}R$ and $\eta\in \mathcal{S}(\R^d)$, we have
\begin{align*}
\sum_{ z\in\Z^d,|z|\geq10R^\e}\eta\Bigl(\frac{y}{R}-z\Bigr) \lesssim
_N\int_{|x|>5R^\e}\Bigl(1+\left|\frac{y}{R}-x\right|\Bigr)^{-N}dx
\lesssim_N R^{-\e(N-d)}.
\end{align*}
Choosing $N\thickapprox\frac{100d}{\e}+d$, we obtain
$$\Big|\sum\limits_{ |z|\geq10R^\e}\eta\left(\frac{y}{R}-z\right)\Big|\lesssim_\e R^{-100d}.$$
By Cauchy-Schwarz and \eqref{osc-sol}, we have
$$
|J_{1,1}|\lesssim_\e R^{-100
d}\|\mathcal{F}_{y\rightarrow\xi}K(x,\cdot,t_j)\|_{L^2_\xi}\|Tf(\cdot,t_j)\|_{L^2_{y}}\lesssim_\e
R^{-100 d}\|f\|_2.
$$
Thus $J_{1,1}$ is negligible.
\vskip0.4cm

$\bullet$\textbf{ Estimation of $J_{1,2}$}.
Notice first that the phase
function of $K(x,y;t_j)$ reads
$$\Phi(x,\,y,\,\xi,\,t_j)\mathrel{\mathop=^{\rm def}}(x'-y)\cdot\xi+(x_{d+1}-t_j)|\xi|^\a.$$
Thus the critical points of $\Phi(x,\,y,\,\xi,\,t_j)$ occurs only when
$$y=x'+(x_{d+1}-t_j)\a|\xi|^{\a-2}\xi.$$
Since $R^\e\gg \a2^{\a+1}>1$ and
$$|x'|\leq 2R,\,0<x_{d+1}-t_j<R,\,
|y|\geq \a 2^{\a+2}R,$$
we have by triangle inequality
$$
\big|\nabla_\xi\Phi\big|\geq|y|-|x'|-\a 2^{\a-1}|x_{d+1}-t_j|\geq\, |y|-\a
2^{\a+1}R\geq R.
$$
Using integration by parts, we may estimate $K(x,y;t_j)$ in $J_{1,2}$  by
$$
|K(x,\,y;\,t_j)|\lesssim_\a\psi(x')(|y|-\a 2^{\a+1}R)^{-100d}.
$$
As a consequence, we have $$|J_{1,2}|\lesssim_\a \psi(x') R^{-99d}\|f\|_2.$$
Hence this term is also negligible.

\subsection{ The estimation of $J_2$.}
Now, let us start to estimate $J_2$.
Rewrite $J_2$ as
\begin{align*}
  J_2=\psi(x')\int e^{i[ x'\cdot\xi+(x_{d+1}-t_j)|\xi|^\a]}\chi(\xi)\^{f_j}(\xi)d\xi
     =\psi(x')Tf_j(x',x_{d+1}-t_j),
\end{align*}
where we have adopted the following notion
\begin{equation}\label{def fj}
f_j(y)\mathrel{\mathop=^{\rm def}}\sum_{ y^\sigma\in\mathfrak{X},\,
|y^\sg|<10R^{1+\e}}\eta\Bigl(\frac{y-y^\sg}{R}\Bigr)Tf(y,t_j).
\end{equation}
This gives the first term on the right side of \eqref{lem 1 '}.
In the sequel, it suffices to show the $f_j$'s, defined by \eqref{def fj}, satisfy \eqref{lem 1''}.
\vskip0.2cm

To do this, we perform some reductions first.
Let $\{\xi^{(k)}\}_k$ be a family of maximal $R^{1-\a}-$separated points
of $\mathcal{I}_{\frac{1}{R}}$ and
cover $\mathcal{I}_{\frac{1}R}$ with essentially disjoint balls
$B(\xi^{(k)},\,R^{1-\a})$. This covering admits a partition
of unit as follows
$$\sum_k \varphi_k(\xi)=1,$$
where $\varphi_k$ is a smooth function supported in the ball $B(\xi^{(k)},R^{1-\a})$.
On account of this, we may write $f=\sum\limits_k f_{(k)}$ and
$f_j=\sum\limits_kf_{j,\,(k)}$ for $j\in\{0,\ldots,\,[R^{\a-1}]\}$,
where
$$\^f_{(k)}(\xi)=\^f(\xi)\varphi_k(\xi),\, \text{and}\q\^f_{j,\,(k)}(\xi)=\^f_{j}(\xi)\varphi_k(\xi),$$
are all supported in $B(\xi^{(k)},\,R^{1-\a})$.
By Plancherel's theorem and almost orthogonality, it
suffices to find some $c_d>0$ such that
\begin{equation}\label{2.3}
    \sum^{[R^{\a-1}]}_{j=0}\|f_{j,\,(k)}\|^2_2\leq C_\e R^{\e c_d}\|f_{(k)}\|^2_2,
\end{equation}
with $C_\e>0$ independent of $k$.\vskip0.2cm

Without loss of generality, we only deal with the case when $k=0$ and suppress the
subscript~$k$ in $f_{(k)}$ and $f_{j,\,(k)}$ for brevity. As a result,
we may assume $\text{supp }\^f\subset B(\xi^{(0)},\,R^{1-\a})$
in the following argument and
normalize $\|f\|_2=1$. By Cauchy-Schwarz, we have
\begin{align*}
    |f_j(y)|^2
    \lesssim R^{2d\e}\sum_{y^\sg\in\mathfrak{X},
    |y^\sg|<10R^{1+\e}}\eta^2\Bigl(\frac{y-y^\sg}{R}\Bigr)\bigl|Tf(y,t_j)\bigr|^2.
\end{align*}
Integrating both sides with respect to $y$ and summing up over $j$, we obtain
\begin{align}\label{2.4}
   \sum_j\|f_j\|^2_2\lesssim R^{2 d\e} \sum_j\sum_{y^\sg\in\mathfrak{X},|y^\sg|<10R^{1+\e}}\int\eta^2\Bigl(\frac{y-y^\sg}{R}\Bigr)|Tf(y,t_j)|^2dy.
\end{align}
Invoking the definition of $Tf(y,t_j)$, we can write
$$
Tf(y,t_j) =I_1+I_2,
$$
where
\begin{align*}
 I_1\mathrel{\mathop=^{\rm def}}&\int_{\Omega^c_{y,\,j}}\mathfrak{K}(y,z;t_j)f(z)\,dz,\\
 I_2 \mathrel{\mathop=^{\rm def}}&\int_{\Omega_{y,\,j}}\mathfrak{K}(y,z;t_j)f(z)\,dz,\\
 \mathfrak{K}(y,z;t_j)\mathrel{\mathop=^{\rm def}}&\int e^{i[(y-z)\cdot\xi+t_j|\xi|^\a]}\chi(\xi)\,d\xi,\\
 \Omega_{y,\,j}\mathrel{\mathop=^{\rm def}}&\bigl\{z\in\R^d\mid\;\bigl|z-y-\a t_j|\xi^{(0)}|^{\a-2}\xi^{(0)}\bigr|<\a 2^{\a+2} R\bigr\}.
\end{align*}
Thus, it suffices to evaluate the contribution of $I_1$ and $I_2$ to \eqref{2.4}.
\vskip0.2cm

 $\bullet$\textbf{ The contribution of $I_1$}.
\vskip0.2cm
Since  $|\xi-\xi^{(0)}|\leq R^{1-\a}$ and $|t_j|\leq R^\a$, we have
\begin{align*}
&\Big|\nabla_\xi\big[(y-z)\cdot\xi+t_j|\xi|^\a\big]\Big|\\
\geq&\big|z-(y+\a t_j|\xi^{(0)}|^{\a-2}\xi^{(0)})\big|-\a2^{\a}R\\
\geq&\a 2^{\a+1} R.\end{align*}
This allows us to use integration by parts to evaluate
$\mathfrak{K}(y,z;t_j)$
$$\big|\mathfrak{K}(y,z;t_j)\big|\lesssim_N\big|z-(y+\a t_j|\xi^{(0)}|^{\a-2}\xi^{(0)})\big|^{-N}.$$
Choosing $N$ large enough, we see the contribution of $I_1$ to
\eqref{2.4} is bounded by
\begin{align*}
&R^{2\varepsilon
d}\sum_j\sum_{y^\sg\in\mathfrak{X},|y^\sigma|\leq10R^{1+\varepsilon}}\int_{\R^d}\eta\Big(\frac{y-y^\sigma}R\Big)^2\big|I_1\big|^2\,dy
\\
\lesssim_N\,&R^{3\varepsilon
d}\sum_j\sup_{y^\sigma}\int_{\R^d}\eta\Big(\frac{y-y^\sigma}R\Big)^2\int_{\Omega_{y,j}^c}\big|z-(y+\a
t_j|\xi^{(0)}|^{\a-2}\xi^{(0)})\big|^{-2N}dzdy
\\
\lesssim_N\,&R^{3\varepsilon d-2N+2d+\a-1} \lesssim_\e\,R^{-200d}.
\end{align*}

$\bullet$\textbf{ The contribution of $I_2$}.
\vskip0.2cm
We use Poisson's summation formula with respect to $z-$variable in $I_2$ to get
$$\sum_j\sum_{y^\sg\in\mathfrak{X},|y^\sigma|\leq10R^{1+\varepsilon}}
\int_{\R^d}\eta\Big(\frac{y-y^\sigma}R\Big)^2\biggl|\int_{\Omega_{y,j}}\mathfrak{K}(y,z;t_j)f(z)dz\biggr|^2dy\leq
L_1+L_2, $$ where
\begin{align*}
L_1=&\sum_j\sum_{\substack{y^\sigma\in\mathfrak{X}\\|y^\sigma|\leq10R^{1+\varepsilon}}}\int_{\R^d}\eta\Big(\frac{y-y^\sigma}R\Big)^2
    \biggl|\sum_{\substack{z_0\in\mathfrak{X}\\ z_0\notin\mathfrak{A}(y)}}\int_{\Omega_{y,j}}\eta\Big(\frac{z-z_0}{R}\Big)\mathfrak{K}(y,z;t_j)f(z)dz\biggr|^2dy\\
L_2=&\sum_j\sum_{\substack{y^\sigma\in \mathfrak{X}\\|y^\sigma|\leq10R^{1+\varepsilon}}}
    \int_{\R^d}\eta\Big(\frac{y-y^\sigma}R\Big)^2\bigg|\sum_{\substack{z_0\in\mathfrak{X}\\ z_0\in\mathfrak{A}(y)}}
    \int_{\Omega_{y,j}}\eta\Big(\frac{z-z_0}{R}\Big)\mathfrak{K}(y,z;t_j)f(z)dz\bigg|^2dy,\\
&\mathfrak{A}(y)\mathrel{\mathop=^{\rm
   def}}\bigl\{z_0\in\mathfrak{X}\mid |z_0-(y+\a\, t_j|\xi^{(0)}|^{\a-2}\xi^{(0)})|\leq
   10R^{1+\varepsilon}\bigr\}.
\end{align*}

Now, we show $L_1$ is also negligible. In fact, since
$$|z_0-(y+\a
t_j|\xi^{(0)}|^{\a-2}\xi^{(0)})|> 10R^{1+\varepsilon},$$ and
$$\big|z-(y+\a
t_j|\xi^{(0)}|^{\a-2}\xi^{(0)})\big|<\a2^{\a+2}R,$$
we have by letting $R$ sufficiently large such that $R^\e\gg \a2^{\a+2}$
$$|z-z_0|\geq\frac12|z_0-(y+\a t_j|\xi^{(0)}|^{\a-2}\xi^{(0)})|.$$
Under the above conditions, we obtain
\begin{align*}
\sum_{z_0\in\mathfrak{X}\setminus\mathfrak{A}(y)}\eta
\Big(\frac{z-z_0}{R}\Big)
&\lesssim_N R^N\sum_{
z_0\in\mathfrak{X}\setminus\mathfrak{A}(y)}|z_0-(y+\a
t_j|\xi^{(0)}|^{\a-2}\xi^{(0)})|^{-N}\\
&\lesssim_N R^{-(N-d)(1+\e)}.
\end{align*}
Choosing
$$N\approx\frac1{2\e}(200 d+\a+1),$$
and using Cauchy-Schwarz as before, we get
\begin{align*}
L_1\lesssim_N &R^{-2\e
N+2d(1+\e)}\sum_j\sum_{|y^\sigma|\leq10R^{1+\varepsilon}}\int_{\R^d}\eta\Big(\frac{y-y^\sigma}R\Big)^2\bigg(
\int_{\Omega_{y,j}}\big|\mathfrak{K}(y,z;t_j)f(z)\big|dz\bigg)^2\,dy
\\
\lesssim_N&R^{-2\e N+5d+\a-1}\lesssim_\e R^{-100d}.
\end{align*}
Thus the contribution of $L_1$ is negligible.
\vskip 0.2cm
Next, we turn to the evaluation of $L_2$. This term contains the nontrivial contribution of \eqref{2.4}.
First, applying Cauchy-Schwarz's inequality to the summation with respect to $z_0$, we have
$$L_2\lesssim\,R^{\e d}(H_1+H_2),$$
where $H_\gamma$ is defined as follows, for $\gamma=1,2$
\begin{align*}
H_\gamma=\sum_j\sum_{|y^\sigma|\leq10R^{1+\varepsilon}}\sum_{z_0\in\mathfrak{A}(y)}\int_{\R^d}\eta\Big(\frac{y-y^\sigma}R\Big)^2
\bigg|\int_{\Omega^\gamma_{y,j,z_0}}\eta\Big(\frac{z-z_0}{R}\Big)\mathfrak{K}(y,z;t_j)f(z)dz\bigg|^2dy
\end{align*}
with $$\Omega^1_{y,j,z_0}=\Omega_{y,j}\cap
B(z_0,\,R^{1+\varepsilon}),\;\Omega^2_{y,j,z_0}=\Omega_{y,j}\setminus
B(z_0,\,R^{1+\varepsilon}).$$

\vskip0.2cm

$\bullet$ {\bf The evaluation of $H_1$:} In view of
Fubini's theorem and
$$|t_j-t_{j+1}|=R,\,\q j=0,\ldots,[R^{\a-1}]-1,$$ there are at most
$R^\varepsilon$ many  $\Omega_{y,j}$'s to intersect with
$B(z_0,R^{1+\varepsilon})$. Denote by $\chi_1(z)$  the
characteristic function of $B(z_0,\,R^{1+\e})$. Applying
Plancherel's theroem to $H_1$ yields
\begin{align*}
H_1
\lesssim &R^{\varepsilon
}\max_j\sum_{z_0\in\mathfrak{X}}\int_{\R^d}\bigg|\iint
e^{i[(y-z)\cdot\xi+t_j|\xi|^\a]}\chi_{\Omega_{y,j}}
\eta\Big(\frac{z-z_0}{R}\Big)\chi_1(z) f(z)dzd\xi\bigg|^2dy
\\
\lesssim &R^{\varepsilon
}\max_j\sum_{z_0\in\mathfrak{X}}\int_{\R^d}\Bigg|\int
e^{i(y\cdot\xi+t_j|\xi|^\a)}
\mathcal{F}_{z\rightarrow\xi}\bigg(\chi_{\Omega_{y,j}}\eta\Big(\frac{\cdot-z_0}{R}\Big)\chi_1
f(\cdot)\bigg)(\xi)d\xi\Bigg|^2dy
\\
\lesssim &R^{\varepsilon
}\sum_{z_0\in\mathcal{X}}\int_{\R^d}\Bigl[\eta\Bigl(\frac{z-z_0}{R}\Bigr)\chi_1(z)
|f(z)|\Bigr]^2dz \lesssim R^{\varepsilon }.
\end{align*}

$\bullet$ {\bf The evaluation of $H_2$:} Since $\eta(x)\in \mathcal{S}(\R^d)$, we have
\begin{align*}
H_2
\lesssim_N &R^{2\varepsilon
d}\sup_{y^\sigma}\sum_j\int_{\R^d}\eta\Big(\frac{y-y^\sigma}R\Big)\sup_{z_0}\bigg|\int_{|z-z_0|\geq
R^{1+\varepsilon}}\Big|\frac{z-z_0}{R}\Big|^{-N}
|f(z)|dz\bigg|^2dy\\
\lesssim_N&R^{-\varepsilon (N-d)+d+\a-1}\lesssim_\e
R^{\varepsilon d}.
\end{align*}
Collecting all the estimations on $I_1$, $L_1$ and $H_1$, we get
eventually \eqref{lem 1''}, and this completes the proof of Lemma
\ref{lem 1}.

\begin{rem}
After finishing this work, we are informed by  Professor Sanghyuk Lee that Lemma \ref{lem 1} can be deduced without losing $R^\e$ by adapting the argument  for the temporal localization Lemma 2.1 in \cite{ref ChoLeeVargas}. We have however decided to include our method since it exhibits different techniques which are interesting in its own right.
\end{rem}

\section{Proof of Lemma \ref{key_decomp}}\label{sect:key decomp}
\noindent

This section is devoted to the proof of Lemma \ref{key_decomp},
which is divided into three parts.
First, we establish an auxiliary result Lemma \ref{lem tri-est}.
Second, we deduce an inductive formula with respect to different scales by
exploiting the self-similarity of Lemma \ref{lem tri-est}.
Finally, we  iterate this inductive formula to get Lemma \ref{key_decomp}.

\subsection{An auxiliary lemma}
Let us begin with an outline of the main steps.
First, we partition the support of
$\^{f}$ into the union of $\frac{1}{K}\times\frac{1}{K}-$ squares with
$K\ll R$. Then, we rewrite $Tf(x)$ as the superposition of the
solutions of the linear Schr\"{o}dinger equation, where each of the initial datum is
frequency-localized in one of these squares. This oscillatory
integral $Tf(x)$ can be transferred into an exponential sum, where the fluctuations of the
coefficients on every box $Q_{a,K}$ of size
$K^{1-\e}\times K^{1-\e}\times K^{1-\e}$ are so slight that they can be viewed essentially
as a constant on each of such boxes.

Next, we partition
$B(0,R)\times[0,\,R]\subset\R^3$ into the union of disjoint $Q_{a,K}$
and estimate the exponential sum on each $Q_{a,K}$. In doing so, we
encounter three possibilities.
For the first one, we will have the transversality condition so that
the multilinear restriction theorem in \cite{ref BCT} can be applied.
To handle the case when the transversality fails, we consider the other two
possibilities. For this part, we use more information from geometric structures
along with Cordorba's square functions \cite{ref Co}. Now, let us turn to details.\vskip0.2cm

Partition\ $\mathcal{I}$ into the union of disjoint
$\frac{1}{K}\times\frac{1}{K}$ squares $\Omega_\nu$, centered at $\xi_\nu$
$$\mathcal{I}\subset\bigcup\limits_{\nu}\Omega_\nu.$$
Then, we rewrite $Tf(x)$ into an exponential sum
$$
Tf(x)=\sum_\nu T_{\nu}f(x)e^{i\phi(x,\xi_\nu)},
$$
where $\phi(x,\xi)=x_1\xi_1+x_2\xi_2+x_3|\xi|^\a$ and
\begin{align*}
T_{\nu}f(x)=\int_{\Omega_\nu}
e^{i[\phi(x,\,\xi)-\phi(x,\,\xi_\nu)]}\^{f}(\xi)\,d\xi.
\end{align*}

$\bullet$ {\bf The local constant property of $T_\nu f(x)$.}
\vskip0.2cm

From a
direct computation, we see $\widehat{T_\nu f}(y)$ is
supported in the following set
$$\left\{y\in\R^3\;\Big|\; y=(y_1,y_2,y_3),\,|y_j|\leq \frac1K, j=1,2,3\right\}.$$
If we take a smooth radial function
$\hat{\eta}(\omega)$ such that $\hat{\eta}(\omega)=1$ for
$|\omega|<2$ and $\hat{\eta}(\omega)=0$ for
$|\omega|>4,\omega\in\R^3$,
then $$ \widehat{T_\nu
f}(\omega)=\widehat{T_\nu f}(\omega)\widehat{\eta_K}(\omega), $$
for $\eta_K(x)=K^{-3}\eta(K^{-1}x)$.
Consequently, $T_\nu f=T_\nu f*\eta_K$.
\vskip0.2cm

Let $Q_a=Q_{a,K}$ be a $K^{1-\e}\times K^{1-\e}\times K^{1-\e}$ box centered
at $a\in K^{1-\e}\Z^3$.
Then, we have
$$B(0,R)\times[0,\,R]\subset\bigcup_a Q_a,$$
where the union is taken over all $a$ such that $Q_a\cap(B(0,R)\times[0,\,R])\neq\emptyset$.
Denote
by $\chi_a=\chi_{Q_a}(x)$ the characteristic function of~$Q_a$.
Restricting $x\in Q_{a}$ and making change of variables
$x=\tilde{x}+a$ with $\tilde{x}\in Q_{0,K}$, we have
\begin{align*}
|T_\nu f(x)|=&\Bigl|\int \tau_{-a}(T_\nu f)(z)\eta_K(\tilde{x}-z)dz\Bigr|\\
\leq&\int \bigl|\tau_{-a}(T_\nu f)(z)\bigr|\sup_{\tilde{x}\in
Q_{0,K}}|\eta_K(\tilde{x}-z)|dz \mathrel{\mathop=^{\rm
def}}c_{a,\nu}.
\end{align*}
Thus, we associate to any $Q_a$ a sequence of $\{c_{a,\nu}\}_{\nu}$.
\vskip0.2cm

Noting that
$$|\nabla_x\phi(x,\xi)-\nabla_x\phi(x,\xi_\nu)|\leq \frac{1}{K},\;\forall\,\xi\in \Omega_\nu,$$
one easily derives the following uniform estimate
$$
|T_\nu f(x')-T_\nu f(x'')|\leq \frac{\|f_\nu\|_2}{K},\;\forall\, |x'-x''|<K.
$$
In particular, $|T_\nu f(x)|$ deviates from $c_{a,\nu}$ by only $\frac{\|f_\nu\|_2}{K}$ whenever $x\in Q_a$.
\vskip0.2cm

The local constant trick can also be regarded as an extension of the Shannon-Nyquist sampling
theorem in \cite{Tao:time-frequency harmonic analysis}.
From
\[T_\nu f(x)=\int T_\nu f(x-Ky)\eta(y)dy,\]
and that $\int \eta=1$, we see that $T_\nu f(x)$ is essentially an average
on the ball $B(x,K)$ up to Schwartz tails. From the uncertainty principle,
$|T_\nu f(x)|$ is essentially constant on the boxes $Q_{a,K}$.
We refer to \cite{Tao:time-frequency harmonic analysis} for standard exposition on this issue.
Thus, whenever $x$ belongs to $Q_a$, we may regard $T_\nu f(x)$ as $c_{\nu,a}$.
\vskip0.2cm

In the preparations, we set $K^{1-\e}$ to be the side length of $Q_{a,K}$ where $\e$ is necessary for a technical reason
so that the local constant property holds. However, to simplify the notations, we will suppress this small $\e$ in the following context. Keeping this in mind, we next classify these $Q_a$'s into three categories.

\vskip0.2cm

$\bullet$ {\bf The classification of $\{Q_a\}_a$.}
\vskip0.2cm
Let $\mathcal{A}$  consist of all  $a$ associated to the boxes $Q_a$ as above.
We will write $\mathcal{A}$ into the union of $\mathcal{A}_j$ for $j=1,2,3$
with $\mathcal{A}_j\subset\mathcal{A}$ defined as follows. \vskip0.2cm

Let $c^*_a=\max\limits_\nu c_{a,\,\nu}$ and $\xi_{\nu^*_a}$ be the center of the square
$\Omega_{\nu^*}$ associated to $c^*_a$.
We define $\mathcal{A}_1\subset\mathcal{A}$ to be such that $a\in
\mathcal{A}_1$ if and only if there exist
$\nu_1,\,\nu_2,\,\nu_3\in\{1,\ldots,\sim K^2\}$ with the property that
$$
\min\{c_{a,\,\nu_1},\,c_{a,\,\nu_2},\,c_{a,\,\nu_3}\}>K^{-4} c^*_a,
$$
and $\xi_{\nu_1},\,\xi_{\nu_2},\,\xi_{\nu_3}$ are non-collinear in
the sense of
\begin{align}\label{non-col}
|\xi_{\nu_1}-\xi_{\nu_2}|\geq|\xi_{\nu_1}-\xi_{\nu_3}|\geq
\text{dist}(\xi_{\nu_3},\,\ell\,(\nu_1,\,\nu_2))>10^3\frac{\a2^\a}K,
\end{align}
where $\ell(\nu_1,\,\nu_2)$ is the straight line through $\xi_{\nu_1}$,
$\xi_{\nu_2}$. \vskip0.2cm

Next, we take $1\ll K_1\ll K\ll R$ and define
$\mathcal{A}_2\subset\mathcal{A}$ such that $a\in\mathcal{A}_2$ if
and only if the following statement holds
\begin{equation}\label{5.2-}
|\xi_{\nu}-\xi_{\nu^*_a}|>4/K_1\Longrightarrow c_{a,\,\nu}\leq K^{-4}c^*_a.
\end{equation}

Let
$ \mathcal{A}_3=(\mathcal{A}\setminus\mathcal{A}_1)\cap(\mathcal{A}\setminus\mathcal{A}_2)$.
Then
$
\displaystyle\mathcal{A}\subset\mathcal{A}_1\cup\mathcal{A}_2\cup\mathcal{A}_3.
$\vskip0.2cm

We claim that if $a\in \mathcal{A}_3$, then there exists a $\nu^{**}_a$
such that $c_{a,\nu^{**}_a}>K^{-4}c^*_a$ and
\begin{equation}\label{5.2*}
|\xi_{\nu^{**}_a}-\xi_{\nu^*_a}|>\frac{4}{K_1},
\end{equation}
moreover, the following statement is valid
\begin{equation}\label{A 3}
\text{dist}(\xi_\nu,\,\ell(\nu^*,\nu^{**}))>10^3\frac{\a2^\a K_1}K\,
\Longrightarrow\, c_{a,\nu}\leq K^{-4} c^*_a.
\end{equation}
Thus, $\{Q_a\}_a$ is classified according to $a$ being inside of
$\mathcal{A}_1$, $\mathcal{A}_2$ or $\mathcal{A}_3$.
\vskip0.2cm

Since the first part of the claim is clear from the definition of $\mathcal{A}_3$,
it suffices to show \eqref{A 3}. We will draw this by
contradiction. Suppose there is an $a\in\mathcal{A}_3$ for which
\eqref{A 3} fails, then there is a $\nu$ such that
$c_{a,\nu}>K^{-4}c^*_a$ but
\begin{equation}\label{A 3'}
\text{dist}(\xi_\nu,\,\ell(\nu^*_a,\nu^{**}_a))>10^3\frac{\a2^\a K_1}K.
\end{equation}
We claim \eqref{A 3'} implies  that  $\xi_\nu,\,\xi_{\nu^*_a}$ and
$\xi_{\nu^{**}_a}$ satisfy \eqref{non-col}. Hence $a\in
\mathcal{A}_1$, which contradicts to
$\mathcal{A}_1\cap\mathcal{A}_3=\emptyset$ and we conclude \eqref{A
3}.
To prove the claim, we consider the following two alternatives due to
symmetry,
\begin{itemize}
\item¡¡\text{case 1:}
$\min\{|\xi_\nu-\xi_{\nu^{**}_a}|,\,|\xi_\nu-\xi_{\nu^*_a}|\}\leq
|\xi_{\nu^*_a}-\xi_{\nu^{**}_a}|$;
\end{itemize}

\begin{itemize}
\item¡¡\text{case 2:}
$\min\{|\xi_\nu-\xi_{\nu^{**}_a}|,\,|\xi_\nu-\xi_{\nu^*_a}|\}>
|\xi_{\nu^*_a}-\xi_{\nu^{**}_a}|$.
\end{itemize}
For case 1, assuming
$|\xi_\nu-\xi_{\nu^{*}_a}|\leq|\xi_\nu-\xi_{\nu^{**}_a}|$
without loss of generality,
we conclude \eqref{non-col}  immediately from
\begin{align*}
|\xi_{\nu^*_a}-\xi_{\nu^{**}_a}|\geq&|\xi_\nu-\xi_{\nu^*_a}|
\geq\text{dist}\bigl(\xi_\nu,\,\ell(\nu^*_a,\nu^{**}_a)\bigr)\\>&10^3\a 2^\a
\frac{K_1}K>10^3 \frac{\a2^\a}K.
\end{align*}

\begin{figure}[ht]
\begin{center}
$$\ecriture{\includegraphics[width=6cm]{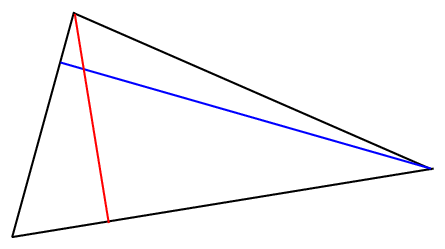}}
{\aat{-2}{-2}{$\xi_{\nu^{**}_a}$}\aat{9}{10}{$h$}\aat{10}{28}{$\xi_{\nu^*_a}$}\aat{26}{11}{$H$}\aat{52}{5}{$\xi_{\nu}$}}$$
\end{center}
\caption{The triangle  $\triangle\,
\xi_{\nu}\xi_{\nu^*_a}\xi_{\nu^{**}_a}$}
\end{figure}
For case 2, we may assume
$|\xi_{\nu^*_a}-\xi_{\nu^{**}_a}|<|\xi_\nu-\xi_{\nu^{*}_a}|\leq|\xi_\nu-\xi_{\nu^{**}_a}|$
by symmetry and denote $\triangle\,
\xi_{\nu}\xi_{\nu^*_a}\xi_{\nu^{**}_a}$ as the triangle formed by $\xi_\nu$, $\xi_{\nu^*_a}$ and $\xi_{\nu^{**}_a}$
(Figure 2), with
$$H=\text{dist}\bigl(\xi_{\nu},\,\ell(\nu^*_a,\nu^{**}_a)\bigr),$$
$$h=\text{dist}\bigl(\xi_{\nu^*_a},\,\ell(\nu,\nu^{**}_a)\bigr).$$
Considering the measure of
$\triangle\,\xi_{\nu}\xi_{\nu^*_a}\xi_{\nu^{**}_a}$, we get from \eqref{5.2*}
$$
h
=\frac{|\xi_{\nu^*_a}-\xi_{\nu^{**}_a}|}{|\xi_{\nu^{**}_a}-\xi_{\nu}|}H\geq
\frac{4}{K_1}\cdot\frac{10^3K_1}{K}\a2^\a\cdot\frac{1}{4}>10^3 \frac{\a2^\a}K,
$$
where we have used $|\xi_{\nu^{**}_a}-\xi_{\nu}|\leq 4$ and hence
\eqref{non-col} follows.

\vskip0.2cm
 $\bullet$ {\bf An auxiliary lemma.}
\vskip0.2cm

The the following lemma, which we are ready to prove,
exhibits once more the spirit of Bourgain-Guth's method, namely the
failure of non-coplanar interactions implies small Fourier supports with possible
addition separateness structures.
\begin{lem}\label{lem tri-est}
Let $B(0,R)\times[0,R]\subset \cup_a Q_a$ be as before. We have on each $Q_a$
\begin{align}\label{3.0}
|Tf(x)|\lesssim &K^8\max_{\substack{\nu_1,\,\nu_2,\,\nu_2
\\
\text{non-collinear}}}\Bigl(\prod^3_{j=1}|T_{\nu_j}f|\Bigr)^{1/3}(x)
\\
\label{3.2}&+\max_\mu\Bigl|\int_{\tilde{\Omega}_\mu}e^{i\phi(x,\xi)}\hat{f}(\xi)d\xi\Bigr|
\\
\label{3.1}&+K^2_1\max_{\substack{\mathcal{L}_1,\mathcal{L}_2\subset\mathcal{L}_a\\\text{dist}(\mathcal{L}_1,\mathcal{L}_2)>1/K_1}}
\prod^2_{j=1}\Bigl|\sum_{\Omega_\nu\subset\mathcal{L}_j}e^{i\phi(x,\xi_\nu)}T_{\nu}f\Bigr|^{\frac{1}{2}}(x)
\\
\label{3.3}&+K^3_1 K^{-\frac12}\bigl(\sum_{\Omega_\nu\subset\mathcal{L}_a}|T_\nu f(x)|^2\bigr)^{\frac12},
\end{align}
where $\tilde{\Omega}_\mu$ is a
$\frac{1}{K_1}\times\frac{1}{K_1}$ square centered at $\xi_{\mu}\in
\mathcal{I}$. $\mathcal{L}_a$ is a $(\a 2^\a
10^3K_1/K)$-neighborhood of the line $\ell(\nu^*_a,\nu^{**}_a):=\ell^*_a$.
The two separated portions $\mathcal{L}_1,\,\mathcal{L}_2$ of
$\mathcal{L}_a$ are obtained by intersecting $\mathcal{L}_a$ with some
$\tilde{\Omega}_{\mu_1}$ and $\tilde{\Omega}_{\mu_2}$ respectively.
\end{lem}

\begin{proof}
For $x\in Q_a$, we  estimate $|Tf(x)|$ in different ways according to
$a\in \mathcal{A}_1$, $a\in \mathcal{A}_2$ and $a\in \mathcal{A}_3$.
If $a\in \mathcal{A}_1$, we have  from
the definition of $\mathcal{A}_1$
\begin{align*}
    |Tf(x)|^3&\leq (K^2 c^*_a)^3
     \leq K^6 K^{12}
    c_{a,\,\nu_1}c_{a,\,\nu_2}c_{a,\,\nu_3}\\
    &\leq K^{18}
    \sum_{\substack{\nu_1,\,\nu_2,\,\nu_3\\\text{non-collinear}}}\prod^3_{j=1}|T_{\nu_j}f|(x)\\
    &\leq K^{24}\max_{\substack{\nu_1,\,\nu_2,\,\nu_3\\\text{non-collinear}}}\prod^3_{j=1}|T_{\nu_j}f|(x).
\end{align*}
If $a\in \mathcal{A}_2$, we use the statement \eqref{5.2-} for $\mathcal{A}_2$ to estimate $Tf(x)$.
For brevity, we define $\|\xi\|=\max\{|\xi_1|,\,|\xi_2|\},\forall\xi\in
\R^2$ and let
$$
\Omega_{*}=\bigcup\Bigl\{\Omega_\nu\,\Big| \,\|\xi_\nu-\xi_{\nu^*_a}\|\leq
10/K_1\Bigr\}.
$$
Then $\Omega_* $ is a $\frac{10}{K_1}\times\frac{10}{K_1}-$square with center $\xi_{\nu^*_a}$.
We may write
\begin{align*}
    | Tf(x)|&\leq\Big|
    \sum_{\|\xi_\nu-\xi_{\nu^*_a}\|\leq 10/K_1}e^{i\phi(x,\,\xi_\nu)}T_\nu f(x)\Bigr|
    +\Bigl|\sum_{\|\xi_\nu-\xi_{\nu^*_a}\|>10/K_1}e^{i\phi(x,\,\xi_\nu)}T_\nu f(x)\Big|
    \\
    &\leq
    \Big|\sum_{\|\xi_\nu-\xi_{\nu^*_a}\|\leq
    10/K_1}\int_{\Omega_\nu}e^{i\phi(x,\xi)}\hat{f}(\xi)d\xi\Big|+\sum_{|\xi_\nu-\xi_{\nu^*_a}|>4/K_1}\Big|T_\nu f(x)\Big|
    \\
    &\leq 100\max_{\tilde{\Omega}_\mu\subset\Omega_{\nu^*_a}}\Big|\int_{\tilde{\Omega}_\mu}e^{i\phi(x,\xi)}\hat{f}(\xi)d\xi\Big|+K^2
    K^{-4}c^*_a\\
    &\lesssim\eqref{3.2}+\eqref{3.3}.
\end{align*}
If  $a\in \mathcal{A}_3$, we write
$$Tf(x)=\mathfrak{D}_1+\mathfrak{D}_2,$$
where
$$
\mathfrak{D}_j=\chi_a(x)\int_{D_j}e^{i\phi(x,\xi)}\^{f}(\xi)d\xi,\,j=1,2,$$
with
$$
D_1=\bigl\{\xi\,\mid\,\text{dist}(\xi,\ell^*_a)>\a 2^\a10^4K_1/K\bigr\},\;
D_2=\bigl\{\xi\,\mid\,\text{dist}(\xi,\ell^*_a)\leq\a 2^\a 10^4K_1/K\bigr\}. $$
From \eqref{A 3}, we have
\begin{align*}
   |\mathfrak{D}_1|&\lesssim
   \sum_{\xi_\nu:\text{dist}(\xi_\nu,\ell^*_a)>\a 2^\a10^3K_1/K}\bigl|e^{i\phi(x,\xi_\nu)}T_{\nu}f(x)\bigr|
   \leq K^2 K^{-4} c^*_a\leq\eqref{3.3}.
\end{align*}
To evaluate  $\mathfrak{D}_2$, we assume without loss of generality
$$
\text{supp}\,\^{f}\subset\Bigl\{\xi\in\R^2\;\Big|\; \text{dist}(\xi,\,\ell^*_a)\leq \a
2^\a 10^4\frac{K_1}K\Bigr\}.
$$
Let $\{\tilde\Omega_\mu\}_\mu$ be a family of disjoint $\frac{1}{K_1}\times\frac{1}{K_1}-$squares
such that $\mathcal{I}\subset\cup_{\mu}\tilde\Omega_\mu$.
For any $x\in Q_a$, we define
$$
c_\mu(x)=\int_{\tilde{\Omega}_\mu}e^{i\phi(x,\xi)}\hat{f}(\xi)d\xi.
$$
Let
$$
\mathcal{H}_a\mathrel{\mathop=^{\rm def}}
\Bigl\{x\in Q_{a}\;\Big|\;\bigl|Tf(x)\bigr|\leq10^8\max_\mu\bigl|c_\mu(x)\bigr|\Bigr\}.
$$
Then
$$
|\mathfrak{D}_2|\leq|Tf(x)|\chi_{\mathcal{H}_a}(x)+|Tf(x)|\chi_{Q_a\setminus\mathcal{H}_a}(x),
$$
where the first term is bounded by \eqref{3.2}.
To handle the second term, we observe that $x\in Q_{a}\setminus \mathcal{H}_a$ implies
\begin{equation}\label{all-small}
|c_\mu(x)|\leq 10^{-8}|Tf(x)|,\q\forall\mu.
\end{equation}
Consider the following set
$$
\mathcal{J}(x)=\biggl\{\mu\;\Big|\;\frac{10^{-1}}{K^2_1}|Tf(x)|\leq
|c_{\mu}(x)|\leq 10^{-8}|Tf(x)|,\,x\in
Q_a\setminus\mathcal{H}_a\biggr\}.
$$
We have
$\#\mathcal{J}(x)\geq10^7,\forall\,x\in
Q_a\setminus\mathcal{H}_a.$
Indeed, suppose$\#\mathcal{J}(x_0)<10^7$ for some $x_0$, then
\begin{equation}\label{claim-1}
|Tf(x_0)|\leq\sum_{\mu\in\mathcal{J}(x_0)}|c_\mu(x_0)|+\sum_{\mu\not\in\mathcal{J}(x_0)}|c_\mu(x_0)|.
\end{equation}
Because of \eqref{all-small}, we can bound the right side of \eqref{claim-1} with
$$
10^7\times 10^{-8}|Tf(x_0)|+K^2_1\cdot \frac{10^{-1}}{K^2_1}|Tf(x_0)|<\frac15|Tf(x_0)|,
$$
which is impossible. Note that the centers of $\{\tilde\Omega_\mu\}_\mu$ are $\frac{1}{K_1}-$separated, we can
choose $\mu_1,\,\mu_2\in\mathcal{J}(x)$ ( $x-$dependent) such that
$$\text{dist}(\tilde{\Omega}_{\mu_1},\,\tilde{\Omega}_{\mu_2})\geq10^4/K_1,$$
and
$$
|Tf(x)|\leq 10
K^2_1\min\bigl\{|c_{\mu_1}(x)|,\,|c_{\mu_2}(x)|\bigr\},x\in Q_a\setminus \mathcal H_a.
$$
It follows that
\begin{equation}\label{claim-2}
|Tf(x)|\chi_{Q_a\setminus\mathcal{H}_a}(x)\leq 10K^2_1\prod^2_{j=1}\Bigl|\int_{\tilde\Omega_{\mu_j}}e^{i\phi(x,\xi)}\^f(\xi)d\xi\Bigr|^{\frac12},
\end{equation}
where $\mu_1$ and $\mu_2$ might depend on $x$.
Now in view of \eqref{5.2-}, we estimate further
\begin{equation}\label{claim-3}
|c_{\mu_j}(x)|\leq\Bigl|\sum_{\substack{\Omega_\nu\subset\tilde\Omega_{\mu_j}\\
\text{dist}(\xi_\nu,\ell^*_a)\leq\a2^\a10^3K_1/K}}\int_{\Omega_\nu}e^{i\phi(x,\xi)}\^f(\xi)d\xi\Bigr|+K^2\cdot K^{-4}c^*_a.
\end{equation}
Thus
\begin{align}
\nonumber&|Tf(x)|\chi_{Q_a\setminus\mathcal{H}_a}(x)\\
\label{tri-terms}\leq& 10K^2_1\Bigl[\prod^2_{j=1}\Bigl|\sum_{\Omega_\nu\subset\mathcal{L}_j}e^{i\phi(x,\xi_\nu)}T_{\nu}f\Bigr|^{\frac12}(x)
+\Bigl(2K^{-2}c^*_a\sum_{\Omega_\nu\subset\mathcal{L}_a}|T_\nu f(x)|\Bigr)^{\frac12}+K^{-2}c^*_a\Bigr],
\end{align}
where $\mathcal{L}_j=\tilde\Omega_{\mu_j}\cap\mathcal{L}_a$,
with $\mathcal{L}_a$ an $\a2^\a10^3K_1/K-$neighborhood of $\ell^*_a$.
Note that
$$c^*_a\leq\bigl(\sum_{\Omega_\nu\subset\mathcal{L}_a}|T_\nu f(x)|^2\bigr)^{\frac12},$$ and
$$\sum_{\Omega_\nu\subset\mathcal{L}}|T_\nu f(x)|\lesssim_\a K_1K^{\frac12}\bigl(\sum_{\Omega_\nu\subset\mathcal{L}_a}|T_\nu f(x)|^2\bigr)^{\frac12},$$
the last two terms in \eqref{tri-terms} is bounded by
$$
K^3_1 K^{-\frac12}\bigl(\sum_{\Omega_\nu\subset\mathcal{L}_a}|T_\nu f(x)|^2\bigr)^{\frac12}.
$$
Therefore, we have
\begin{align*}
|\mathfrak{D}_2|\lesssim\eqref{3.2}+\eqref{3.1}+\eqref{3.3}.
\end{align*}
This completes the proof of Lemma \ref{lem tri-est}.
\end{proof}

\subsection{A self-similar iterative formula}

The formula we deduce in this part is the engine for the
iterating process in order for the fractional order Bourgain-Guth's inequality.
Let $\Omega_\t$ be a $\d-$square centered at $\xi_\t$,
\begin{align*}
&\Omega^\t_\nu\mathrel{\mathop=^{\rm def}}\Bigl\{\xi\in\R^2\,\Big|\,\|\xi-(\xi_\t+\d\xi_\nu)\|<\frac\d K\Bigr\},\;\\
&\tilde\Omega^\t_\mu\mathrel{\mathop=^{\rm def}}\Bigl\{\xi\in\R^2\,\Big|\,\|\xi-(\xi_\t+\d\xi_\mu)\|<\frac\d {K_1}\Bigr\}.
\end{align*}
Denote by $\^f_{\t,\nu}=\^f\cdot\chi_{\Omega^\t_\nu}$.
We prove the following iterative formula from scale $\d$ to $\frac\d K$
for all $a\in \mathcal A$
\begin{align}
\label{scaled} |T_\t f|(x)\lesssim&
K^8\max_{\substack{\nu_1,\,\nu_2,\,\nu_3\\\text{non-collinear}}}\prod^3_{j=1}\left|Tf_{\t,\,\nu_j}\right|^{\frac{1}{3}}(x)
\\
\label{scaled_2}
&+\psi_\t(x)\Big(\sum_{\Omega^\t_{\nu}\subset\mathcal{L}_\t}
|Tf_{\t,\,\nu}|^2\Big)^{\frac{1}2}(x)
\\
\label{scaled_3}
&+\max_{\mu}\Bigl|\int_{\tilde\Omega^\t_\mu\cap\Omega_\t}e^{i[\xi\cdot
x'+x_3|\xi|^\a]}\hat{f}(\xi)d\xi\Bigr|,
\end{align}
for all $x\in\T^*_a$,
where $$
\mathcal{T}^*_a\mathrel{\mathop=^{\rm def}}\Bigl\{x\in\R^3\,\Big|\,
|x_1-a_1|<\frac K\d,\,|x_2-a_2|<\frac{K}{\d}, \,
|x_3-a_3|<\frac{K}{\d^\a}\Bigr\},
$$
and $\mathcal{L}_\t$ is a $\d/K-$neighborhood of a line segment. $\psi_\t$ is
a function satisfying
\begin{equation}\label{averagescaled}
\left(\frac{1}{|B|}\int_{B}\psi_\t^4(x)dx\right)^{1/4}\lesssim
K^{2\a}_1,
\end{equation}
where $B$ is a $K\mathcal{C}^*_\t-$box centered at $a$.
\vskip0.2cm

To deduce \eqref{scaled}-\eqref{averagescaled}, we need the following estimate,
which is standard square functions going back to Cordoba \cite{ref Co}.
The crucial $L^4$-estimate is used in \cite{ref Bourgain Guth} to tackle the worst
scenario by exploiting the separateness of the line segments in  which the frequencies locate.
This part of frequencies corresponds to the terms of main contributions.
\begin{lem}For any $a\in\mathcal{A}$ and all $x\in Q_{a}$, we have
\begin{equation}\label{dual1}
    \Bigl(\frac{1}{|Q_{a}|}\int_{Q_{a}}\bigl|\eqref{3.1}+\eqref{3.3}\bigr|^4dx\Bigr)^{1/4}\lesssim
    K^{2\a}_1\Bigl(\sum_{\Omega_\nu\subset\mathcal{L}}|T_{\nu}f(x)|^2\Bigr)^{\frac{1}{2}},
\end{equation}
where the implicit constant is independent of $a$.
\end{lem}

\begin{rem}
This observation is crucial for the iteration in the next subsection.
To prove \eqref{dual1}, we rely heavily on the $10^4/K_1-$separateness of the segments $\mathcal{L}_1,\,\mathcal{L}_2$.
Since we are dealing with the fractional order symbol, the proof is more intricate than that in \cite{ref Bourgain Guth},
where the algebraic structure simplifies the proof significantly.
\end{rem}

\begin{proof} We need to estimate the $L^4-$averagement of
\eqref{3.1} and \eqref{3.3} over $Q_{a}$. Since on every $Q_a$, $|T_\nu f(x)|$
can be viewed as a constant, we immediately get
$$
\Bigl(\frac{1}{|Q_a|}\int_{Q_a}\eqref{3.3}^4dx\Bigr)^{\frac14}
\lesssim K^3_1\Bigl(\sum_{\Omega_\nu\subset\mathcal{L}}|T_{\nu}f(x)|^2\Bigr)^{\frac{1}{2}}.
$$
Next, we estimate \eqref{3.1}. First, we have
\begin{align}\nonumber
&\int_{Q_{a}}\max_{\substack{\mathcal{L}_1,\mathcal{L}_2\subset\mathcal{L}_a\\\text{dist}(\mathcal{L}_1,\mathcal{L}_2)\geq\frac{10^4}{K_1}}}\prod^2_{j=1}
\Bigl|\sum_{\Omega_\nu\subset\mathcal{L}_j}e^{i\phi(x,\xi_\nu)}T_\nu f(x)\Bigr|^2(x)dx\\
\label{dual'}\leq&\sum_{\substack{\mathcal{L}_1,\mathcal{L}_2\subset\mathcal{L}_a\\\text{dist}(\mathcal{L}_1,\mathcal{L}_2)\geq\frac{10^4}{K_1}}}
\int_{Q_{a}\cap\mathfrak{S}_{\mathcal{L}_1,\mathcal{L}_2}}\prod^2_{j=1}
\Bigl|\sum_{\Omega_\nu\subset\mathcal{L}_j}e^{i\phi(x,\xi_\nu)}T_\nu f(x)\Bigr|^2(x)dx
\end{align}
where
$$
 \mathfrak{S}_{\mathcal{L}_1,\,\mathcal{L}_2}\mathrel{\mathop=^{\rm def}}\biggl\{x\,\Big|\; \eqref{3.1}\leq
 2\,\Bigl|\sum_{\Omega_\nu\subset\mathcal{L}_1}e^{i\phi(x,\xi_\nu)}T_{\nu}f\Bigr|^{\frac{1}{2}}
 \Bigl|\sum_{\Omega_\nu\subset\mathcal{L}_2}e^{i\phi(x,\xi_\nu)}T_{\nu}f\Bigr|^{\frac{1}{2}}(x)\biggr\},
$$
and the summation in \eqref{dual'} is taken over all the pairs of $10^4/K_1-$separated subsegments of $\mathcal{L}_a$.
Since there are at most $K^2_1$ pairs of such subsegments,
it suffices to estimate each term in the summation. Take a Schwartz function $\rho\geq 0$ such that $\^\rho$ is compactly supported and
$\rho_a(x):=\rho\big(\frac{x-a}{K}\big)=1$ for all $x\in Q_a$
\begin{align}\label{bilinear1}
&\int_{Q_{a}\cap\mathfrak{S}_{\mathcal{L}_1,\,\mathcal{L}_2}}\prod^2_{j=1}
\Bigl|\sum_{\Omega_\nu\subset\mathcal{L}_j}e^{i\phi(x,\xi_\nu)}T_\nu f(x)\Bigr|^2dx\\
\leq&\int_{\R^3}\prod^2_{j=1}
\Bigl|\rho_a(x)\sum_{\Omega_\nu\subset\mathcal{L}_j}e^{i\phi(x,\xi_\nu)}T_\nu f(x)\Bigr|^2dx
\\
\label{bilinear2}\leq&\sum_{\substack{\Omega_{\nu_1},\,\Omega_{\nu_2}\subset\mathcal{L}_1
\\
\nonumber\Omega_{\nu'_1},\,\Omega_{\nu'_2}\subset\mathcal{L}_2}}
\Bigl|\int_{\R^3}\big[\rho_aT_{\nu_1}f\cdot\overline{\rho_aT_{\nu_2}f\cdot
\rho_aT_{\nu'_1}f}\cdot\rho_aT_{\nu'_2}f\big](x)
e^{i\Psi(x,\,\xi_{\nu_1},\,\xi_{\nu_2},\,\xi_{\nu'_1},\,\xi_{\nu'_2})}dx\Bigr|,
\end{align}
where
$$
 \Psi(x,\,\xi_{\nu_1},\,\xi_{\nu_2},\,\xi_{\nu'_1},\,\xi_{\nu'_2})=
 \phi(x,\xi_{\nu_1})-\phi(x,\xi_{\nu_2})-\phi(x,\xi_{\nu'_1})+\phi(x,\xi_{\nu'_2}).
$$
Considering the support of the function
$$\widehat{\rho_a}*\widehat{T_{\nu_1}f}*\widehat{\rho_a}*
\widehat{\overline{T_{\nu_2}f}}*\widehat{\rho_a}*\widehat{\overline{
T_{\nu'_1}f}}*\widehat{\rho_a}*\widehat{T_{\nu'_2}f},$$
we may restrict the summation  to
those quadruples $\nu_1,\,\nu_2,\,\nu'_1,\,\nu'_2$ such that
\begin{equation}
\begin{cases}
\label{quad1}|\xi_{\nu_1}-\xi_{\nu_2}-\xi_{\nu'_1}+\xi_{\nu'_2}|\lesssim\frac{1}{K},
\\
\big||\xi_{\nu_1}|^\a-|\xi_{\nu_2}|^\a-|\xi_{\nu'_1}|^\a+|\xi_{\nu'_2}|^\a\big|\lesssim\frac{1}{K}.
\end{cases}
\end{equation}
Denote by $\ell^*_a=\R v+b$ where $v,\, b\in\R^2$ with $|v|=1$, $|b|\leq
4$ and $v\perp b$. Since $\xi_{\nu_1},\,\xi_{\nu_2},\,\xi_{\nu'_1}$
and $\xi_{\nu'_2}$ are lying in a $\frac{\a 2^\a10^3K_1}{K}-$ neighborhood
of $\ell^*_a$, there are $t_1,\,t_2,\,t'_1$ and $t'_2$ such that for $j=1,2$, we have
$$|\xi_{\nu_j}-t_jv-b|\leq \a2^{\a}10^3\frac{K_1}{K},\;|\xi_{\nu'_j}-(t'_j)v-b|\leq \a 2^\a 10^3\frac{K_1}{K}.$$
In view of \eqref{quad1} and the $10^4/K_1-$separateness of $\mathcal{L}_1$ and $\mathcal{L}_2$, we have
\begin{equation}\label{cond 1}
|t_1-t_2|<\frac{2}{K_1},\,|t'_1-t'_2|<\frac2
{K_1},\;|t_1-t'_1|>\frac{10^4}{K_1},
\end{equation}
\begin{equation}\label{cond 2}
|t_1-t_2-t'_1+t'_2|\lesssim \frac {K_1} K,
\end{equation}
\begin{equation}\label{cond 3}
|\varphi(t_1)-\varphi(t_2)-\varphi(t'_1)+\varphi(t'_2)|\lesssim\frac
{K_1} K,
\end{equation}
where $\varphi(t):=(t^2+|b|^2)^{\frac{\a}{2}}$. \vskip0.2cm
 We claim \eqref{cond
1}, \eqref{cond 2} and \eqref{cond 3} imply
\begin{equation}\label{orth goal}
|t_1-t_2|\lesssim \frac{K^\a_1}{K},\q|t'_1-t'_2|\lesssim\frac{K^\a_1}K.
\end{equation}
As a consequence, we may deduce that
$$
|\xi_{\nu_1}-\xi_{\nu_2}|\lesssim \frac{K^\a_1}{K},\,|\xi_{\nu'_1}-\xi_{\nu'_2}|\lesssim \frac{K^\a_1}{K},
$$
which implies
$$
\eqref{bilinear1}\lesssim K^{4\a}_1\sum_{\Omega_{\nu}\subset\mathcal{L}_1,\,
\Omega_{\nu'}\subset\mathcal{L}_2}\int_{Q_{a}}|Tf_{\nu}|^2|Tf_{\nu'}|^2(x)dx.
$$
Summing up \eqref{bilinear1} with respect to pairs of $(\mathcal{L}_1,\,\mathcal{L}_2)$,
we obtain from the local constant property fulfilled by the $|T_\nu f|(x)$'s
\begin{align*}
\int_{Q_{a}}|\eqref{3.1}|^4dx
\lesssim & K^{4\a+2}_1\int_{Q_{a}}\Bigl(\sum_{\Omega_\nu\subset\mathcal{L}}|Tf_{\nu}(x)|^2\Bigr)^2dx\\
\lesssim &K^{4\a+2}_1|Q_{a}|\Bigl(\sum_{\Omega_\nu\subset\mathcal{L}}|Tf_{\nu}(x)|^2\Bigr)^2,\q\forall\,x\in Q_{a}.
\end{align*}
 This proves \eqref{dual1}.
\vskip0.2cm

To prove the claim \eqref{orth goal}, we need to consider the following two possibilities.
$$ t_1<t_2<t'_1<t'_2\,,\;\text{and}\q t_1<t_2<t'_2<t'_1.$$
In view of \eqref{cond 2}, the second possibility implies
$$|t_1-t_2|+|t'_1-t'_2|\lesssim\frac{K_1}{K}.$$
Then \eqref{orth goal} follows immediately.
It is thus sufficient to handle the first possibility.
To do this, we consider the following three cases.
\vskip0.2cm $\bullet$ $t_1<t_2<0<t'_1<t'_2$. For
this case, we have from \eqref{cond 3}
\begin{equation}\label{sum}
\bigl|\varphi(t_1)-\varphi(t_2)\bigr|+\bigl|\varphi(t'_1)-\varphi(t'_2)\bigr|\lesssim\frac{K_1}{K}.
\end{equation}
By triangle inequality and \eqref{cond 1}, we have either
$t'_1>\frac{10^3}{K_1}$ or $t_2<-\frac{10^3}{K_1}$.
We only handle the case when $t'_1>\frac{10^3}{K_1}$ since the other one is exactly the same.
We apply mean value theorem to $\varphi(t)$ to get a $t_*$ with
$t'_1<t_*<t'_2$ such that
\begin{align*}
\frac{K_1}{K}\gtrsim\bigl|\varphi(t'_1)-\varphi(t'_2)\bigr| \gtrsim
\a|t_*|^{\a-1}|t'_1-t'_2| \gtrsim  \frac{\a}{K^{\a-1}_1}|t'_1-t'_2|.
\end{align*}
Hence, we have
$$|t'_1-t'_2|\lesssim_\a
\frac{K^{\a}_1}{K}.$$
By triangle inequality again and \eqref{cond 2}, we obtain
$$|t_1-t_2|\lesssim_\a\frac{K^{\a}_1}{K}.$$

\vskip0.2cm $\bullet$ $t_1<0<t_2<t'_1<t'_2$.
First,  we always have $t'>\frac{10^3}{K_1}$ in this case.
If $t_1<-t_2$, we also get \eqref{sum},
so \eqref{orth goal} follows immediately as above.
For $-t_2\leq t_1$, we have from \eqref{cond 3}
\begin{equation}\label{minus}
\Bigl|\bigl|\varphi(t_1)-\varphi(t_2)\bigr|-\bigl|\varphi(t'_1)-\varphi(t'_2)\bigr|\Bigr|\lesssim\frac{K_1}{K}.
\end{equation}
Suppose $|t_1-t_2|\gg \frac{K^\a_1}{K}$,
then $|t'_1-t'_2|\gg \frac{K^\a_1}{K}$ by \eqref{cond 1}.
Moreover, we have
\begin{equation}\label{equivalent}
\frac12|t'_1-t'_2|\leq|t_1-t_2|\leq2|t'_1-t'_2|.
\end{equation}
By mean value theorem, we have for some $t_*\in(t_1,t_2)$, $t'_*\in (t'_1,t'_2)$
with $t'_*>\max\{\frac{10^3}{K_1},\,t_*+\frac{10^3}{K_1}\}$ such that
\begin{align*}
\frac{K_1}{K} &\gtrsim
\Bigl|\bigl(\varphi(t'_2)-\varphi(t'_1)\bigr)-\bigl(\varphi(t_2)-\varphi(t_1)\bigr)\Bigr|\\
&\gg\a\frac{K^\a_1}{K} \Big(|t'_*|^{\a-1}-|t_*|^{\a-1}\Big) \gg \a
\frac{K_1}K,
\end{align*}
which is impossible since $\a>1$ and $K_1\gg 1$.

\vskip0.2cm $\bullet$ $0<t_1<t_2<t'_1<t'_2$. This can be reduced to the above two cases, and
we complete the proof of \eqref{dual1}.
\end{proof}

Here, another observation made by Bourgain and Guth in \cite{ref Bourgain Guth} is that
as a consequence of \eqref{dual1}, for $x\in Q_{a}$ one define $\psi(x)$  by writing
$$
\eqref{3.1}+\eqref{3.3}=\psi(x)\Bigl(\sum_{\Omega_\nu\subset\mathcal{L}}|Tf_\nu(x)|^2\Bigr)^{\frac{1}{2}}.
$$
Clearly $\psi$ is nonnegative such that
\begin{equation}\label{2.12'}
\Bigl(\frac{1}{|Q_{a}|}\int_{Q_{a}}\psi^4(x)dx\Bigr)^{1/4}\lesssim
K^{2\a}_1.
\end{equation}
To see this is possible, one only need to be aware of the local constant property
enjoyed by $T_\nu f(x)$'s so that $\psi(x)$ can be defined on each ball of radius $K$
due to \eqref{dual1}. Then we glue all the pieces of $\psi(x)$ on the balls together.
By an average argument and the local constancy property of $\eqref{3.1}+\eqref{3.3}$ and functions $T_\nu f(x)$
on each $Q_{a,K}$ box, we may assume $\psi(x)$ is constant on the unit cubes centered at lattices.
\\

\begin{rem}
By writing \eqref{3.1} + \eqref{3.3} into a product of an appropriate $\psi$ and a square function,
we may iterate this part step by step in the subsequent context to generate the items having transverality
structures corresponding to all the dyadic scales. This is one of the brilliant ideas due to Bourgain and Guth,
which is also applied by Bourgain in \cite{ref Bourgain3} and \cite{ref Bourgain Guth} as a substitution
of Wolff's induction on scale technique.
\end{rem}

Substituting $\eqref{3.1}+\eqref{3.3}$ in Lemma \ref{lem tri-est} for
$$\psi(x)\Bigl(\sum_{\Omega_\nu\subset\mathcal{L}}|T
f_\nu(x)|^2\Bigr)^{\frac{1}{2}},$$
we obtain
\begin{align}\label{new-tri-est}
|Tf(x)|\lesssim &K^8\max_{\substack{\nu_1,\,\nu_2,\,\nu_2
\\
\text{non-collinear}}}\Bigl(\prod^3_{j=1}|T_{\nu_j}f|\Bigr)^{1/3}(x)
\\
\label{new-tri-est2}&+\psi(x)\Bigl(\sum_{\Omega_\nu\subset\mathcal{L}}|T
f_\nu(x)|^2\Bigr)^{\frac{1}{2}}
\\
\label{new-tri-est3}&+\max_\mu\Bigl|\int_{\tilde{\Omega}_\mu}e^{i\phi(x,\xi)}\hat{f}(\xi)d\xi\Bigr|.
\end{align}
Now, we are ready to prove \eqref{scaled}-\eqref{averagescaled}.
Observe that $Tf(x)$ is controlled in terms of
$Tf_\nu$'s with $\^f_\nu$  supported in a square of size
$\frac1K$, whereas $\^f$ is supported in a region of size $1$. Thus
it is natural to scale each $\^f_\nu$ to be a function $\^g_\nu$
such that supp $\^g_\nu$ is of size $1$. After applying
\eqref{new-tri-est}-\eqref{new-tri-est3} to each
$Tg_\nu$, we  rescale the estimates on $Tg_\nu$ back to the original size $\frac1K$. More
generally, this process can be carried out  with $Tf_\t$ in place of $Tf(x)$
on the left side of \eqref{new-tri-est}, where $\^f_\t$ is supported in
a square of size $\d$.
\begin{proof}[Proof of \eqref{scaled}-\eqref{averagescaled}]
Let
$\^{f}_\t=\^{f}|_{\Omega_\t}$ and
$
x\in\mathcal{T}^*_a=a+Q^\d_{0,K},
$
with $$Q^\d_{0,K}=\{x\,\mid \,(\d x_1,\d x_2,\d^\a x_3)\in Q_{0,K}\}.$$
Making change of variables
\begin{equation}
\label{change-variable} x'=a'+\tilde{x}'/\d,\;
x_3=a_3+\tilde{x}_3/\d^\a,\;\xi=\xi_\t+\d\eta,
\end{equation}
where $\xi_\t$ the center of $\Omega_\t$, and $\tilde{x}\in Q_{0,K}$, we have
\begin{align}
    \nonumber\chi_{\mathcal{T}^*_a}(x)T f_\t(x)
    =&\chi_{Q_{0,K}}(\tilde{x})\int_\Omega
    e^{i[\tilde{x}'\cdot\eta+\tilde{x}_3|\eta|^\a]}\widehat{g^{\t,\d}_a}(\eta)d\eta\\
    \label{dot}=&\chi_{Q_{0,K}}(\tilde{x}) T(g^{\t,\d}_a)(\tilde{x}),
\end{align}
with
$$\widehat{g^{\t,\d}_a}(\eta)=e^{i[\d
a'\cdot\eta+a_3\d^\a|\eta|^\a+(a'+\tilde{x}'/\d)\cdot\xi_\t+(a_3+\tilde{x}_3/\d^\a)
(|\xi_\t+\d\eta|^\a-|\d\eta|^\a)]}\d^2\widehat{f}_\t(\xi_\t+\d\eta)\chi_{\Omega}(\eta).$$
Now that $\^{g^{\t,\d}_a}$ is supported in a square of size $1$, we can apply
\eqref{new-tri-est}-\eqref{new-tri-est3} to \eqref{dot} with $\tilde{x}\in Q_{0,K}$ to obtain
\begin{align}\label{doot}
   \Bigl|T(g^{\t,\d}_a)(\tilde{x})\Bigr|
   \lesssim& K^8\max_{\substack{\nu_1,\nu_2,\nu_3\\
   \text{non-collinear}}}\prod^3_{j=1}\bigl|T_{\nu_j}(g^{\t,\d}_a)\bigr|^{\frac13}(\tilde x)\\
   \label{doot'}&+\psi(\tilde{x})\Bigl(\sum_{\Omega_\nu\subset\mathcal{L}}
   |T_{\nu}(g^{\t,\d}_a)|^2\Bigr)^{\frac{1}{2}}(\tilde x)\\
   \label{doot''}&+\max_\mu\bigl|\int_{\widetilde{\Omega}_\mu}e^{i\phi(\tilde{x},\xi)}\^{g^{\t,\d}_a}(\eta)d\eta\bigr|.
\end{align}
Re-scale the $\widehat{g^{\t,\d}_a}$ in \eqref{doot}-\eqref{doot''} back to $\^f_\tau$
by using \eqref{change-variable} and  setting $\eta=(\zeta-\xi_\t)/\d$
\begin{align}
   \nonumber\Bigl| \chi_{Q_{0,K}}(\tilde{x})T_{\nu}(g^{\t,\d}_a)(\tilde{x})\Bigr|
    &=\Bigl|\int_{\|\eta-\xi_{\nu}\|<1/K}e^{i[\tilde{x}'\cdot\eta+\tilde{x}_3|\eta|^\a]}\^{g^{\t,\d}_a}(\eta)d\eta\Bigr|\\
    \nonumber&=\Bigl|\int_{\|\z-(\xi_\t+\d\xi_{\nu})\|<\d/K}e^{i[\z\cdot(a'+\tilde{x}/\d)+(a_3+\frac{\tilde{x}_3}{\d^\a})|\z|^\a]}\hat{f}_\t(\z)d\z\Bigr|\\
    \label{dooot}&=\Bigl|\int_{\|\z-(\xi_\t+\d\xi_{\nu})\|<\d/K}e^{i[\z\cdot x'+x_3|\z|^\a]}\hat{f}_\t(\z)d\z\Bigr|,\;x\in\mathcal{T}^*_a.
\end{align}
From \eqref{dot}, \eqref{doot} and \eqref{dooot}, we get \eqref{scaled}-\eqref{scaled_3} on
$\mathcal{T}^*_a$.

Since
$$
\psi_\t(x)=\psi\bigl(\d(x_1-a_1),\,\d(x_2-a_2),\,\d^\a(x_3-a_3)\bigr),\forall\;x\in\T^*_a,
$$
we obtain \eqref{averagescaled}  from \eqref{2.12'}.
\end{proof}

\subsection{Iteration and the end of the proof}
This part follows intimately the idea of Bourgain and Guth in \cite{ref Bourgain Guth}
however, we provide more explicit calculations during the iteration process so that
this marvelous idea can be reached even for the novice readers. It is hoped that
this robust machine will be upgraded so that further improvements seems possible in this area.
\\

Let $1\ll K_1\ll K\ll R$.
From Lemma \ref{lem
tri-est}, we have
\begin{align}\label{5.1}
    |Tf(x)|\lesssim&\, K^8 \max _{\substack{\Omega_{\t_1},\,\Omega_{\t_2},\,\Omega_{\t_3}:\frac{1}{K}-\text{cubes}\\
    \text{non-collinear}}}\prod^3_{j=1}|T_{\t_j}f|^{\frac{1}{3}}(x)\\
    \label{5.2}&+\psi(x)\Big[\sum_{\substack{\Omega_\t\subset\mathcal{L}\\\Omega_\t:\frac{1}{K}-\text{cubes}}}|T_\t f|^2\Big]^{1/2}(x)\\
   \label{5.2'} &+\max_{\tilde{\Omega}_{\tilde{\t}}:\frac{1}{K_1}-\text{cubes}}|T_{\tilde{\t}}f|(x),
\end{align}
where $\psi$ is approximately constants on unit boxes and obeys
$$\bigg(\frac{1}{|Q_{a}|}\int_{Q_{a}}\psi^4(x)dx\bigg)^\frac14\lesssim K_1^{2\a},$$
for any $Q_{a}$.

Noting that \eqref{5.1} involves a triple product of  $|T_{\t_j} f|^{\frac13}$ with $j=1,2,3$,
we call this term is of type I. For \eqref{5.2}, it is a product of a suitable function $\psi$
and an $\ell^2$ sum of $\{T_\t f\}_\t$, and we call this term is of type II.
The term \eqref{5.2'} is an $\ell^{\infty}$ norm of $\{T_{\widetilde{\t}}f\}_{\widetilde{\t}}$,
and we call it of type III.
In each step of the iteration below, we will encounter plenty
of terms belonging to type I, II and III from the previous step.
These are called \emph{newborn }terms.
We add the newborn  terms of type I to the type I terms of the previous generations
and keep on iterating all the newborn terms of type II and III
to get the next generation of type I, II and III terms. This is the iterating mechanism.

To be more precise, we use
\eqref{scaled} - \eqref{scaled_3} to $T_\t f$ in
\eqref{5.2} with $\d=1/K$ and to $T_{\tilde{\t}} f$ in \eqref{5.2'},
with $\d=1/K_1$. This is exactly the first step of the iteration.
After this, we obtain terms of type I, type II and type III generated by \eqref{5.2} and
\eqref{5.2'}. In each type of the terms, the supports of $\^f_{\t}$'s could be of the
scales like
$$\frac{1}{K^2},\q\frac{1}{KK_1}\q\text{or}\q\frac{1}{K^2_1}.$$
Adding the newborn terms of type I to the previous type I terms,
we repeat the same argument as in the first
step to all the terms of type II and type III to get the second
generation. This process is continued with newborn terms
of type I added to the pervious type I terms until the scale of the support of
$\^f_\t$ in the terms of type II and type III becomes
$\frac{1}{\sqrt{R}}$. Finally, we obtain a collection of type I terms at different scales
and a remainder consisting of type II and III terms at scale $\frac{1}{\sqrt{R}}$,
which is controlled by \eqref{K2}. This yields \eqref{K1} and \eqref{K2}.\vskip0.2cm

Now, we present the explicit computation for the first step.
Applying \eqref{scaled} - \eqref{scaled_3} to
\eqref{5.2} with $\d=1/K$, we have by Minkowski's inequality
\begin{align}\nonumber
      &\label{5.3}\q\psi(x)\Bigl(\sum_{\Omega_\t\subset\mathcal{L}}|T_\t
      f|^2\Bigr)^{1/2}(x)
      \\
     \lesssim & K^8\Bigg[\sum_{\Omega_{\t}\subset\mathcal{L}}
      \Bigg(\max_{\substack{\Omega_{\t}\supset\Omega_{\t^{(1)}_1},\,\Omega_{\t^{(1)}_2},\,\Omega_{\t^{(1)}_3}(x)
      \\\frac{1}{K^2}-\text{squares};\,\text{non-collinear}}}\psi\prod^3_{j=1}|T_{\t^{(1)}_j}f|^{\frac{1}{3}}\Bigg)^2\Bigg]^{\frac{1}{2}}(x)
      \\
      &\label{5.4}\q+\Bigg[\sum_{\Omega_{\t}\subset\mathcal{L}}\sum_{\substack{\mathcal{L}^{(1)}\supset\Omega_{\t^{(1)}}\\
      \frac{1}{K^2}-\text{squares}}}\bigl(\psi\psi_{\t}|T_{\t^{(1)}}f|\bigr)^2\Bigg]^{\frac{1}{2}}(x)
      \\
      &\label{5.5}\q+\Bigg(\sum_{\Omega_{\t}\subset\mathcal{L}}\psi^2\max_{\substack{\tilde{\Omega}_{\tilde{\t}^{(1)}}\subset\Omega_\t
      \\\frac{1}{K_1K}-\text{squares}}}|T_{\tilde{\t}^{(1)}}f|^2\Bigg)^{\frac{1}{2}}(x),
\end{align}
where the superscript in $\t^{(k)}$ denotes the $k-$th step of the iteration.
\vskip0.2cm

Denoting $\psi_{:\t^{(1)}:}=\psi\psi_\t$, we need to verify for any $\e>0$
\begin{equation}\label{L4average}
 \frac{1}{|\C^*_\t|}\int_{\C^*_\t}\psi^4(x)dx\lesssim R^\e,
\end{equation}
and
\begin{equation}\label{L4averageo1st}
 \frac{1}{|\C^*_{\t^{(1)}}|}\int_{\C^*_{\t^{(1)}}}\psi_{:\t^{(1)}:}^4(x)dx\lesssim R^\e.
\end{equation}
To get \eqref{L4average}, we use the boxes $ Q_a$ to subdivide $\C^*_\t$ such that
$$
\C^*_\t\subset\bigcup_a{Q_a}\subset 2\C^*_\t.
$$
Then \eqref{2.12'} gives
\begin{align*}
\frac{1}{|\C^*_\t|}\int_{\C^*_\t}\psi^4(x)dx&\lesssim\frac{1}{|\bigcup_a
Q_{a}|} \int_{\cup Q_{a}}\psi^4(x)dx\\
&\lesssim
\max_{\C^*_\t\subset Q_a\subset2\C^*_\t}\frac{1}{|Q_a|}\int_{Q_a}\psi^4(x)dx
\lesssim K^{8\a}_1\ll R^\e.
\end{align*}
To verify \eqref{L4averageo1st}, we  note that $\C^*_{\t^{(1)}}$ is
a $K^2\times K^2\times K^{2\a}-$box
in the direction of the normal vector of the surface at $\xi_{\t^{(1)}}$.
It follows that $\C^*_{\t^{(1)}}$ can be covered as follows
$$
\C^*_{\t^{(1)}}\subset\bigcup_\t B^*_\t\subset 2\C^*_{\t^{(1)}}
$$
where $B^*_\t$ is a $K\C^*_\t-$box.
Then, we have
\begin{equation}\label{eq5.49}
\frac{1}{|\C^*_{\t^{(1)}}|}\int_{\C^*_{\t^{(1)}}}\psi^4\psi^4_\t(x)dx
\lesssim \max_{B^*_\t\subset2\C^*_{\t^{(1)}}}\frac{1}{|B^*_\t|}\int_{B^*_\t}\psi^4\psi^4_\t(x)dx.
\end{equation}
To estimate the term on the right side, we let $\{\mathcal{B}_\rho\}_\rho$
be a collection of essentially disjoint $\mathcal{C}^*_\t-$boxes  such that
$$B^*_\t\subset\bigcup_\rho\mathcal{B}_\rho\subset 2B^*_\t.$$
Since $\psi_\t$ is approximately constant on each $\mathcal{B}_\rho$, we have
\begin{align}
\nonumber\int_{B^*_\t}\psi^4\psi^4_\t(x)dx&\lesssim\sum_\rho\Big[\int_{\mathcal{B}_\rho}\psi^4(x)dx\Big](\psi_\t|_{\mathcal{B}_\rho})^4,\\
\nonumber\lesssim& K^{8\a}_1\sum_\rho(\psi_\t|_{\mathcal{B}_\rho})^4|\mathcal{B}_\rho|
\lesssim K^{8\a}_1\int_{B^*_\t}\psi^4_\t(x)dx\\
\label{*}\lesssim& K^{16\a}_1|\C^*_\t|\ll R^\e|\C^*_\t|,\end{align}
where we have used \eqref{averagescaled}, \eqref{L4average}
and $\mathcal{B}_\rho$ is a $\mathcal{C}^*_\t-$box.
This along with \eqref{eq5.49} proves \eqref{L4averageo1st}.
\vskip0.2cm

Next, we apply
\eqref{scaled} - \eqref{scaled_3} to \eqref{5.2'}
with $\d=\frac{1}{K_1}$, and obtain
\begin{align}
\q|(\ref{5.2'})| \label{2.22}&\lesssim
\max_{\tilde{\Omega}_\t:\frac{1}{K_1}-\text{squares}}
             K^8\max_{\substack{\tilde{\Omega}_{\tilde{\t}}\supset\Omega_{\tilde{\t}^{(1)}_1},\,\Omega_{\tilde{\t}^{(1)}_2},\,\Omega_{\tilde{\t}^{(1)}_3}\\
             \text{non-collinear};\frac{1}{K_1 K}-\text{squares}}}\prod^3_{j=1}|T_{\tilde{\t}^{(1)}_j}f|^{\frac{1}{3}}\\
\label{2.23}&\q +\max_{\tilde{\Omega}_\t:\frac{1}{K_1}-\text{squares}}
             \Bigg(\sum_{\substack{\Omega_{\tilde{\t}^{(1)}}\subset\mathcal{L}^{(1)}\cap\tilde{\Omega}_{\tilde{\t}}\\\frac{1}{K_1 K}-\text{squares}}}\psi^2_{\tilde{\t}}|T_{\tilde{\t}^{(1)}}f|^2  \Bigg)^{\frac{1}{2}}\\
\label{2.24}&\q +\max_{\tilde{\Omega}_\t:\frac{1}{K_1}-\text{squares}}
\max_{\substack{\tilde{\Omega}_{\t}\supset\tilde{\Omega}_{\t^{(1)}}\\\tilde{\Omega}_{\t^{(1)}}:\frac{1}{K^2_1}-\text{squares}}}|T_{\tilde{\t}^{(1)}}f|.
\end{align}
We also have estimates on the $L^4-$averagement akin to \eqref{L4average} and \eqref{L4averageo1st}.

\begin{rem}
After the first step, if we already have $$K^2\sim
KK_1\sim\sqrt{R}.$$ Then  Lemma \ref{key_decomp} is
proved right after the first iteration. However, this is not the
case since $1\ll K_1\ll K\ll R$. Therefore we have to use the above
argument recursively to conclude Lemma \ref{key_decomp}.
\end{rem}

Noting that the scales for type I terms at the $k-$th generation is $K^{-m_1}K^{-m_2}_1$ with $m_1,m_2\in \Z, \,m_1+m_2=k+1, m_1,m_2\geq0$,
we find the type I terms of the $k-$th stage from
the previous $(k-1)-$th stage is dominated by a $k-$fold sum
$$
\Bigg(\sum_{\Omega_{\t}\subset\mathcal{L}}\sum_{\Omega_{\t^{(1)}}\subset\mathcal{L}^{(1)}}\cdots\sum_{\Omega_{\t^{(k)}}\subset\mathcal{L}^{(k)}}
\max_{\substack{\Omega_{\t^{(k+1)}_j}\subset\Omega_{\t^{(k)}},j=1,2,3\\
\Omega_{\t^{(k+1)}_j}:\frac{1}{K^{k+1}}-\text{squares}\\
\text{non-collinear}}}
\psi^2_{:\t^{(k)}:}\Bigg(\prod^3_{j=1}|T_{\t^{(k+1)}_j}f|\Bigg)^{\frac{2}{3}}\Bigg)^{\frac{1}{2}}+\sum_{\text{mixed scales}},
$$
where the summation over mixed scales represents the cases when $m_2\geq1$.
For brevity, we only write out the case when $m_2=0$ explicitly and the other cases are similar.
Notice that in each fold of the sum, there are at most $KK_1$ terms involved,
the above expression can be controlled by
$$
C(K)\max_{\mathcal{E}_{(\frac{1}{K})^k}}\Bigg[\sum_{\Omega_\t^{(k)}\in
\mathcal{E}_{(\frac{1}{K})^{k(1+\e)}}}\Bigg(\psi_{:\t^{(k)}:}\prod^3_{j=1}|T_{\t^{(k+1)}_j}f|^{\frac{1}{3}}\Bigg)^2\Bigg]^{\frac{1}{2}}+\sum_{\text{mixed scales}},
$$
where we have adopted the notations in Lemma \ref{key_decomp}.
\vskip0.2cm

\begin{figure}[ht]\label{converingrelation}
\begin{center}
$$\ecriture{\includegraphics[width=8cm]{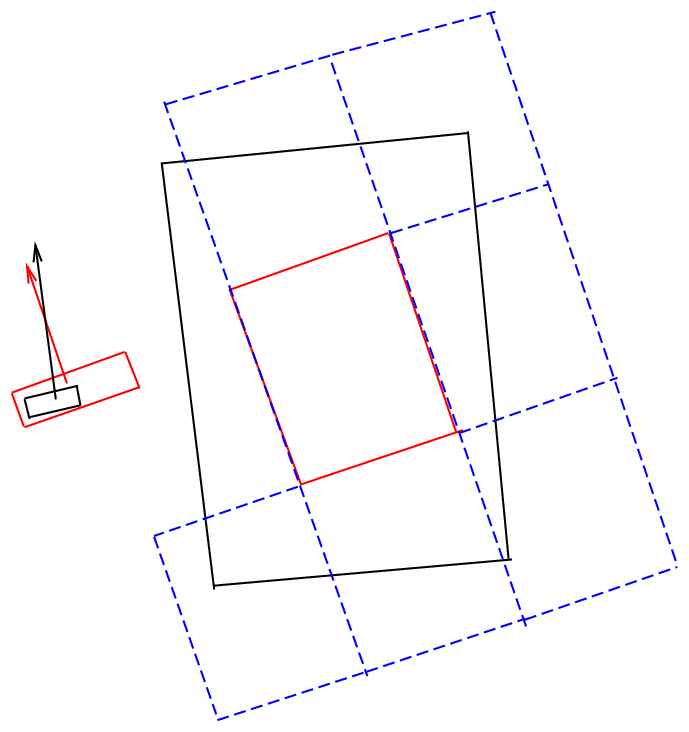}}
{\aat{-4}{23}{$\mathcal{C}_{\t^{(\ell)}}$}\aat{15}{23}{$\mathcal{C}^*_{\t^{(\ell)}}$}
\aat{8}{29}{$\mathcal{C}_{\t^{(\ell-1)}}$}\aat{23}{25}{$K\mathcal{C}^*_{\t^{(\ell-1)}}$}
\aat{39}{15}{$K\mathcal{B}_\rho$}}$$
\end{center}
\caption{The boxes $\mathcal{C}^*_{\t^{(\ell)}}$
and $K\mathcal{C}^*_{\t^{(\ell-1)}}$.}
\end{figure}

Ir remains to show
$$
\frac{1}{|\mathcal{C}^*_{\t^{(k)}}|}\int_{\mathcal{C}^*_{\t^{(k)}}}\psi^4_{:\t^{(k)}:}(x)dx\lesssim R^\e,\q\text{for any }\e>0.
$$
To prove this, we use the induction on $k$.
We observe that the $L^4$
average of $\psi_{:\t^{(1)}:}$ over $\mathcal{C}^{*}_{\t^{(1)}}$ is
bounded by $K^{4\a}_1$ in the first step.
Assuming in the $(\ell-1)-$th stage, we
already have
\begin{equation}\label{lth-induct-average}
\Bigl(\frac{1}{|\mathcal{C}^*_{\t^{(\ell-1)}}|}\int_{\mathcal{C}^*_{\t^{(\ell-1)}}}\psi^4_{:\t^{(\ell-1)}:}\Bigr)^{\frac14}
\lesssim K_1^{2\,\ell\,\a}.
\end{equation}
Since
$\psi_{:\t^{(\ell)}:}=\psi_{:\t^{(\ell-1)}:}\psi_{\t^{(\ell)}}$, we
need to evaluate
\begin{equation}\label{lth-target}
\frac{1}{|\mathcal{C}^*_{\t^{(\ell)}}|}\int_{\mathcal{C}^*_{\t^{(\ell)}}}\psi^4_{:\t^{(\ell-1)}:}\psi^4_{\t^{(\ell)}}(x)dx.
\end{equation}
Because the angle between the two normal vectors of
$\mathcal{C}_{\t^{(\ell-1)}}$ at $\xi_{\t^{(\ell-1)}}$ and
$\mathcal{C}_{\t^{(\ell)}}$ at $\xi_{\t^{(\ell)}}$ is also
controlled by $\frac{1}{K}$( see Figure 3), we  may construct a cover of $\mathcal{C}^*_{\t^{(\ell)}}$ by
$K\mathcal{C}^*_{\t^{(\ell-1)}}-$boxes as follows.
Denote by $\{\mathcal{B}_\rho\}_\rho$ a collection of $\C^*_{\t^{(\ell-1)}}-$boxes such that
(see Figure 3)
$$
\C^*_{\t^{(\ell)}}\subset\bigcup_\rho K\mathcal{B}_\rho\subset2\,\C^*_{\t^{(\ell)}}.
$$
On account of this, we can estimate
\begin{align*}
|\eqref{lth-target}|\lesssim\max_{\mathcal{B}_\rho}\frac{1}{|K\mathcal{B}_\rho|}\int_{K\mathcal{B}_\rho}\psi^4_{:\t^{(\ell-1)}:}\psi^4_{\t^{(\ell)}}(x)dx.
\end{align*}
By the hypothesis of induction, we have
\begin{align*}
\int_{K\mathcal{B}_{\rho}}\psi^4_{:\t^{(\ell-1)}:}\psi^4_{\t^{(\ell)}}(x)dx
&\lesssim\sum_{\rho}(\psi_{\t^{(\ell)}}|_{\mathcal{B}_\rho})^4\int_{\mathcal{B}_{\rho}}\psi^4_{:\t^{(\ell-1)}:}(x)dx\\
&\lesssim K_1^{8\a\ell}\sum_\rho(\psi_{\t^{(\ell)}}|_{\mathcal{B}_\rho})^4|\mathcal{B}_\rho|\\
&\lesssim K_1^{8\a\ell}\int_{K\mathcal{B}_{\rho}}\psi^4_{\t^{(\ell)}}(x)dx \\
&\lesssim K_1^{8\a(\ell+1)}|K\mathcal{B}_{\rho}|,
\end{align*}
where in the last estimate, we used
$$\frac{1}{|K\mathcal{C}^*_{\t^{(\ell-1)}}|}\int_{K\mathcal{C}^*_{\t^{(\ell-1)}}}\psi^4_{\t^{(\ell)}}(x)dx\lesssim K_1^{8\a}.$$
We denote
$$\d=K^{-(\ell+1)},$$
and assume at the $\ell-$th stage
$$\frac{1}{\sqrt{R}}<\d.$$
Noting that $K_1\ll K$ and
$$\ell+1=\log\frac{1}{\d}/\log K,$$
we have
$$
K_1^{8\a(\ell+1)}<R^{\frac{\log K^{4\a}_1}{\log K}}\ll R^\e,\;\forall\;\e>0.
$$
If $\d$ ranges from all the dyadic numbers between $R^{-1/2}$ and $K^{-1}$,
we see the contribution from all the type I terms is bounded by \eqref{K1}.
The contributions from
\eqref{5.5}\eqref{2.23} and \eqref{2.24} to \eqref{K1}
can be evaluated in a similar manner.\\

When the scale arrives at $\frac1{\sqrt{R}}$,
the remainder term is  bounded by \eqref{K2}. Finally, we lose an
$R^\e-$factor by taking maximum in \eqref{K1} and
\eqref{K2} with respect to dyadic $\d\in(R^{-1/2},\,1/K)$.
Thus, we complete the proof of Lemma \ref{key_decomp}.

\subsection*{Acknowledgements}
We would like to appreciate Professor Sanghyuk Lee for indicating to us the recent work \cite{ref ChoLeeVargas}.
The authors were supported by
the NSF of China under grant No.11231006 and 11371059.


\begin{thebibliography}{HD}




\normalsize
\baselineskip=17pt

%
%
%
%
%
%



\bibitem[B1]{ref Bourgain1}
J.~Bourgain.
\emph { A remark on Schr\"{o}dinger operators},
 Isreal J. Math., 77(1992), 1--16.


\bibitem[B2]{ref Bourgain2}
J.~Bourgain,
\emph { Some new estiamtes on osillatory integrals},
Essays on Fourier Analysi in Honor of Elias M. Stein
 Princeton, NJ 1991.
 Princeton Math. Ser., Vol. 42, Princeton University Press,
New Jersey, (1995), 83--112.



\bibitem[B3]{ref Bourgain3}
J.~Bourgain,
\emph { On the Schr\"{o}dinger maximal function in higher dimensions},
Proceedings of the Steklov Institute of Mathematics,
280(2013),  46--60.


\bibitem[BG]{ref Bourgain Guth}
J.~Bourgain and  L.~Guth,
\emph{\em Bounds on oscillatory integral operators based on multilinear estimates},
\newblock {\em Geom. Funct. Anal.}, 21(2011), 1239--1295.

\bibitem[BCT]{ref BCT}
J.~Bennett, T.~Carbery and T.~Tao,
\emph { On the multilinear restriction and Kakeya conjectures},
 Acta Math., 196(2006), 261--302.

\bibitem[Cr]{ref Car}
 L.~Carleson,
\emph {Some analystic problems related to statistical mechanics in
Euclidean Harmonic Analysis},
 Lecture Notes in Math.,779,
 Springer Berlin(1979), 5--45.

\bibitem[CLV]{ref ChoLeeVargas}
C.~Cho, S.~ Lee and A.~Vargas,
\emph{ Problems on pointwise convergence of solutions to the Schr\"{o}dinger equation},
 J. Fourier Anal. Appl. 18 (2012), no. 5, 972--994.


\bibitem[Co]{ref Co}
A.~Cordoba,
\emph { Geometric Fourier analysis},
 Ann. Inst. Fourier , 32(1982), 215--226.

\bibitem[DK]{ref DK}
B.~Dahlberg and C.E.~Kenig.
\emph {A note on the almost everywhere behavior of solutions to the Schr\"{o}dinger equations},
 Harmonic Analysis Lecture Notes in Math., 908,
 Springer, Berlin, (1982), 205--209.

\bibitem[H]{ref Hor}
L.~H\"{o}rmander,
\emph{ Oscillatory integrals and multipliers on $FL^p$},
 Arkiv Math., 11(1973), 1--11.

\bibitem[KPV]{ref KPV}
C.~E.~Kenig, G.~Ponce and L.~Vega,
\emph{ Oscillatory integral and regularity of dispersive equations},
 Indiana University Math. Journal, 40(1991), 33--69.

\bibitem[L]{ref Lee}
S.~Lee.
\emph{  On pointwise convergence of the solutions to Schr\"{o}dinger equations in $\R^2$,}
 Int. Math. Res. Not. (2006), 1--21.

\bibitem[LR]{ref Lee Rogers}
S.~Lee and K.~M.~Rogers,
\emph{  The Schr\"{o}dinger equation along curves and the quantum harmonic oscillator},
 Adv. Math. 229 (2012), no.3, 1359--1379.

\bibitem[M]{ref Miyachi}
A.~Miyachi.
\emph{ On some singular Fourier multipliers,}
 J. of the Faculty of Science. The Univ. of
Tokyo. Sec. IA, 28(1981), 267--315

\bibitem[MVV]{ref M-V-V}
 A.~Moyua, A.~Vargas and L.~Vega.
\emph{ Schr\"{o}dinger maximal function and restriction properties of the Fourer transform},
 IMRN , (1996), 793--815.

\bibitem[R]{Rogers08}
 K. Rogers,
 \emph{ A local smoothing estimate for the
Schr\"odinger equation},
 Advances in Mathematics, 219(2008),
2105--2122.

\bibitem[S]{ref Sj}
P.~Sj\"{o}lin.
\emph{ Regularity of solutions to the Schr\"{o}dinger equation},
 Duke Mathematical Journal, 55(1987), 699-715.

\bibitem[Sh]{ref shao}
S.~Shao.
\emph { On Localization of the Schr\"{o}dinger maximal operator,}
arXiv:1006.2787v1.


\bibitem[T1]{Tao:time-frequency harmonic analysis}
T.~Tao.
\emph {Time-frequency harmonic analysis},
http://www.math.ucla.edu/~tao/254a.1.01w/blurb.html

\bibitem[T2]{ref Tao 1}
T.~Tao. \emph{ A sharp bilinear restriction estiamte for
parabloids},
 Geom. Funct. Anal., 13(2003), 1359--1384.


\bibitem[TV1]{ref T-V1}
T.~Tao and A.~Vargas.
\emph{ A bilinear approach to cone multiplier.I.Restriction estiamtes}.
 GAFA , 10(2000), 185--215.

\bibitem[TV2]{ref T-V2}
 T.~Tao and A.~Vargas.
 \emph{ A bilinear approach to cone multiplier.II. Applications}.
  GAFA , 10(2000), 216--258.


\bibitem[TVV]{ref T-V-V 1998}
 T.~Tao, L.~Vega and A.~Vargas,
 \emph{ A bilinear approach to the restriction and Kakeya conjectures}
 J. Amer. Math. Soci. 11(1998), 967--1000.

\bibitem[V]{ref Vega}
L.~Vega,
\emph{ Schr\"{o}dinger equations: pointwise convergence to the initial data},
 Proceedings of the American
Mathematical Society, 102(1988), 874--878.


\bibitem[W]{ref Wolff}
T.~Wolff,
\emph{A sharp bilinear cone restriction estiamte},
Ann. of Math., (2),153(3)(2001):661--698.







\end{thebibliography}
\end{document}